\documentclass[11pt]{article}
\usepackage[utf8]{inputenc}
\usepackage[T1]{fontenc}
\usepackage{amsmath, amsthm, amsfonts, amssymb}

\usepackage{xcolor}
\usepackage{graphicx} 
\usepackage{multirow, multicol}
\usepackage{booktabs}
\usepackage{colortbl}
\usepackage{subcaption}

\usepackage[onehalfspacing]{setspace}

\usepackage{algorithm, algpseudocode}
\usepackage[skins]{tcolorbox}
\usetikzlibrary{decorations.pathreplacing, positioning, arrows, shapes, arrows.meta, calc, patterns}

\usepackage[margin=1in]{geometry}

\usepackage[final,bookmarksnumbered,breaklinks,colorlinks]{hyperref}
\hypersetup{urlcolor=blue,linkcolor=blue,citecolor=blue}

\usepackage{bm}
\usepackage{physics}
\renewcommand{\var}[1]{{#1}}
\newcommand{\mtr}[1]{\boldsymbol{#1}}
\renewcommand{\Re}[1]{\operatorname{Re}\left(#1\right)}
\renewcommand{\Im}[1]{\operatorname{Im}\left(#1\right)}

\newcommand{\Yt}{\tilde Y}
\newcommand{\Ybt}{\tilde{\boldsymbol{Y}}}
\newcommand{\x}{\boldsymbol{x}}

\newcommand{\bsa}{\boldsymbol{a}}
\newcommand{\bsb}{\boldsymbol{b}}
\newcommand{\bsd}{\boldsymbol{d}}
\newcommand{\bsf}{\boldsymbol{f}}
\newcommand{\bsg}{\boldsymbol{g}}
\newcommand{\bsq}{\boldsymbol{q}}
\newcommand{\bsr}{\boldsymbol{r}}
\newcommand{\bsv}{\boldsymbol{v}}
\newcommand{\bsy}{\boldsymbol{y}}
\newcommand{\bsz}{\boldsymbol{z}}

\newcommand{\bsA}{\bm{A}}
\newcommand{\bsB}{\bm{B}}
\newcommand{\bsE}{\bm{E}}
\newcommand{\bsM}{\bm{M}}
\newcommand{\bsP}{\bm{P}}
\newcommand{\bsQ}{\bm{Q}}
\newcommand{\bsW}{\bm{W}}
\newcommand{\bsX}{\bm{X}}
\newcommand{\bsY}{\bm{Y}}
\newcommand{\bsV}{\bm{V}}

\newcommand{\Wb}{\boldsymbol{W}}

\usepackage{todonotes}

\definecolor{darkgreen}{RGB}{0,160,0}

\usepackage{relsize}

\renewcommand{\t}[1]{#1^{\text{{\relsize{-2}$\hspace{-0.01em}\mathsf{T}$}}}}
\usepackage{authblk}

\usepackage[]{apacite}

\usepackage{mathtools}

\allowdisplaybreaks

\title{Solving Sparse Mixed-Integer Quadratic Problems: Application to the Unit Commitment Problem with Optimal Power Flow}

\author[1,2]{Ignacio Gómez-Casares}
\author[3]{Pietro Belotti}
\author[4,5]{Bissan Ghaddar}
\author[1,2]{Julio González-Díaz}

\affil[1]{Department of Statistics, Mathematical Analysis and Optimization and MODESTYA Research Group, University of Santiago de Compostela, Santiago de Compostela, Spain}
\affil[2]{CITMAga (Galician Center for Mathematical Research and Technology), 15782 Santiago de Compostela, Spain}
\affil[3]{ Department of Electronics, Information and Bioengineering, Politecnico di Milano, Milan, Italy}
\affil[4]{Ivey Business School, Western University, London, Ontario, Canada}
\affil[5]{IE University, Madrid, Spain}

\begin{document}
\maketitle

\begin{abstract}
Mixed-Integer Quadratically Constrained Quadratic Programs (MIQCQP) arise in a variety of applications, particularly in energy, water, and gas systems, where discrete decisions interact with nonconvex quadratic constraints. These problems are computationally challenging due to the combination of combinatorial complexity and nonconvexity, often rendering traditional exact methods ineffective for large-scale instances. We propose a solution framework for sparse MIQCQPs that integrates semidefinite programming relaxations with chordal decomposition techniques to exploit sparsity. By leveraging problem structure, we significantly reduce the size of the semidefinite constraints into smaller, tractable blocks, improving the scalability of the relaxation and the overall branch-and-bound procedure.  We evaluate our framework on the Unit Commitment problem with AC Optimal Power Flow constraints, a practically relevant and highly challenging problem that couples discrete generation decisions with nonlinear AC network physics. Computational results on standard IEEE test cases with up to 118 buses show that our approach provides strong bounds and high-quality solutions, while scaling significantly better than state-of-the-art global optimization solvers.
\end{abstract}




\section{Introduction}
Mixed-Integer Quadratically Constrained Quadratic Programming (MIQCQP) represents a class of optimization problems characterized by the presence of both continuous and discrete decision variables, where the objective function and constraints are quadratic. MIQCQP problems are generally NP-hard due to the combination of nonconvex quadratic terms and integrality restrictions. This combination significantly increases computational complexity, requiring advanced solution techniques. MIQCQPs arise naturally in numerous real-world applications, particularly in energy systems \shortcite{kocuk2017,bingane2019,amorosi2025}, water distribution networks \shortcite{bragalli2012,faria2011}, gas networks \shortcite{li2024,liers2021}, and engineering design \shortcite{misener2009}. 

When solving MIQCQPs, traditional methods, such as standard branch-and-bound, struggle to scale as the problem size increases. Consequently, a substantial body of work has focused on strengthening the continuous relaxations through convexification techniques. One widely used approach is Semidefinite Programming (SDP) relaxation \cite{anjos2011handbook, anstreicher2009semidefinite}, where quadratic terms are lifted into matrix variables and the nonconvex quadratic structure is replaced by a positive semidefinite constraint. This relaxation, often referred to as the Shor relaxation, provides tight bounds but can become computationally expensive for large-scale instances due to the size of the semidefinite matrices. To improve tractability, researchers have proposed Second-Order Cone Programming (SOCP) relaxations \cite{kim2003exact}, which approximate the semidefinite constraint using rotated cone constraints or exploit convex quadratic envelopes, providing a compromise between bound quality and computational efficiency. Another important class of techniques is the Reformulation–Linearization Technique (RLT) \cite{sherali2013reformulation}, which generates linear valid inequalities by multiplying existing constraints and variable bounds and then linearizing the resulting bilinear terms using auxiliary variables. RLT constraints effectively capture variable interactions and significantly strengthen linear relaxations when embedded within branch-and-bound frameworks. 
While recent years have seen progress in convex relaxations and global optimization techniques for quadratic programming, the computational burden remains substantial, especially for large-scale MIQCQPs.


This paper contributes to the literature on SDP relaxations for nonconvex quadratic constraints, yielding Mixed-Integer Semidefinite Programs (MISDP). As mentioned above, MISDPs are computationally demanding due to the potentially large semidefinite matrix variables and the coexistence of discrete and semidefinite constraints. To mitigate this limitation, we exploit the sparsity and structural properties commonly present in real-world instances. In many MIQCQPs, constraints usually involve only a small number of monomials, leading to term sparsity \shortcite{wang2020,Wang:2022}, while interactions among variables are often localized, giving rise to correlative sparsity that can be represented by sparse graphs \shortcite{lasserre2006convergent, waki2006sums}. By exploiting these sparsity patterns, chordal decomposition techniques allow for replacing large semidefinite constraints with blocks of smaller ones defined on the maximal cliques of the sparsity graph. This decomposition substantially reduces the size of the semidefinite matrices, lowering the computational cost of SDP relaxations and improving the scalability of branch-and-bound algorithms that incorporate them.

In this article, we consider MIQCQPs formulated as\footnote{Notation: sets are denoted by capital letters. Following conventions in Power Engineering, which provides the case study for this paper, scalars appear in normal font and may be either lower- or upper-case. Vectors and matrices, also in lower- or upper-case, are denoted in bold.}
%

\begin{subequations}
\label{eq:miqcqp-original}
\begin{spreadlines}{2pt} 
\begin{align}
\textbf{[MIQCQP]} \quad \min_{\x, \bsy} \quad & \t{\x} \bsQ \x + \t{\bsd} \x + \t{\bsr} \bsy \\
\text{s.t.} \quad & \t{\x} \bsA_i \x + \t{\bsa_i} \x + \t{\bsb_i} \bsy \leq c_i, & \forall i \in \mathcal{I} \\
& \t{\bsf_i} \x + \t{\bsg_i} \bsy \leq h_i & \forall i\in \mathcal{I}'\\
& \x \in \mathbb{R}^n, \quad \bsy \in \mathbb Z^m, 
\end{align}
\end{spreadlines}
\end{subequations}
and consisting of the following elements: $\x \in \mathbb{R}^n $ are continuous decision variables, $ \bsy\in \mathbb Z^m$ are integer decision variables, $ \bsQ $ is a matrix defining the quadratic objective function, $\bsd$ and $\bsr$ are linear coefficient vectors, $ \bsA_i $ defines quadratic constraints, $\bsa_i$, $\bsb_i$, and $c_i $ are linear coefficients, and $\bsf_i$, $\bsg_i$, and $h_i$ define additional linear constraints. We assume that the continuous variables $\x$ do not form bilinear terms with the integer variables $\bsy$, resulting in a partially decoupled structure that facilitates the use of conic relaxations and the exploitation of the problem structure. In particular, without mixed bilinear terms of the form $x_iy_j$, the quadratic constraints remain separable from the integrality decisions, allowing convex relaxations such as SOCP or SDP to be applied directly without requiring additional linearization or lifting procedures to handle binary-continuous products. The latter assumption is made without loss of generality for the relevant case in which $\bsy$ are binary variables: any interaction between continuous and binary variables, such as those appearing in mixed bilinear terms, can be reformulated using not only standard big-$M$ techniques but also more sophisticated methods \shortcite{bertsimas2021unified,wei2024convex}. These reformulations preserve the problem structure, while isolating the combinatorial decisions in $\bsy$ from the continuous variables $\x$, resulting in more efficient relaxations. 

The main contribution of this paper is to integrate different methodologies into a unified solution framework for solving large-scale sparse MIQCQPs. More precisely, this is achieved by combining:
\begin{itemize}
\item Conic relaxations: particularly semidefinite relaxations, that preserve tightness while maintaining tractability.
\item Decomposition techniques: which exploit correlative and term sparsity in the quadratic forms to reformulate the problem over smaller overlapping semidefinite matrices.

\item High-quality feasible solution generation: using local solvers to guide the branch-and-bound process effectively and generate high-quality feasible solutions.
\end{itemize}

As a case study, we apply our framework to the Unit Commitment problem with AC Optimal Power Flow constraints, hereafter referred to as UC-OPF, a fundamental energy systems application that yields challenging MIQCQPs with a sparse network structure. We demonstrate that our approach produces high-quality solutions on instances derived from standard IEEE test systems. The framework is designed to scale to large problem instances and efficiently produce high-quality feasible solutions for the UC-OPF and related optimization problems. While the primary focus of this work is on computational performance and practical solution quality rather than on guaranteeing global optimality, the proposed method nevertheless yields an optimality gap, which is typically small for the instances considered in this study. Importantly, we also show that our approach scales significantly better than two state-of-the-art solvers, Gurobi \shortcite{gurobi2025} and BARON \shortcite{baron2025}.

The structure of the paper is organized as follows: Section \ref{sec:literature} reviews some of the main approaches used in the literature to solve MIQCQPs, Section \ref{sec:methodology} describes the solution approach, Section \ref{sec:ucacopf} presents the unit commitment problem with optimal power flow constraints, and Section \ref{sec:results} presents and analyzes the computational results on the UC-OPF. Finally, Section \ref{sec:conclusion} summarizes the main contributions and concludes with future research directions.

\section{Literature Review}\label{sec:literature}

Consider a generic problem instance described by [MIQCQP] above.
%
If $\bsQ$ and all of the $\bsA_i$'s above are positive semidefinite (PSD), then the above problem becomes a {\em convex} MIQCQP, due to the convexity of the continuous relaxation, for which several efficient approaches exist. If there are integer variables (\emph{i.e.}, $m > 0$) then the classical branch-and-bound algorithm can be used, where a lower bound at every branch-and-bound node can be obtained via a solver for convex quadratic optimization, such as an interior-point method, or via Outer Approximation \shortcite{duran1986outer}. Convex MIQCQP problems can also be solved by reformulating them as Mixed-Integer Second-Order Conic Programming (MISOCP) problems \shortcite{vielma2017extended}.

If any of $\bsQ$ or the $\bsA_i$ above are not PSD, then [MIQCQP] is nonconvex. This difficult class of problems is the main target of this article, hence we summarize here a vast amount of work that has been conducted in the past decades for the nonconvex case. A quadratic form $\t{\x}\bsQ\x+\t{\bsd}\x=\sum_{j\in N}\sum_{k\in N} q_{jk}x_j x_k + \sum_{j\in N} d_j x_j$ can be reformulated as the linear expression $\sum_{j\in N}\sum_{k\in N} q_{jk}w_{jk} + \sum_{j\in N} d_j x_j$ by creating extra variables $w_{jk}$ and adding the nonconvex constraints $w_{jk} = x_j x_k$ for all $j\in N$, $k\in N$ where $N:=\{1,2,\ldots{},n\}$. The quadratic form can then be re-written as the linear expression $\mathrm{tr}(\bsQ \bsW) + \t{\bsd}\x$ for convenience, where $\bsW$ is the $n\times n$ matrix defined as $\bsW=(w_{jk})_{j\in N, \, k\in N}$ and $\mathrm{tr}(\bsA  \bsB) = \bsA\bullet \bsB$ is the inner product of two matrices $\bsA$ and $\bsB$ of equal dimensions. Note that $w_{jk} = w_{kj}$.
 
Finding tight lower bounds for [MIQCQP] revolves around the ability to find a tight relaxation of the set linking $\bsW$ and $\x$ variables, which can be written in matricial form: $X = \{(\x,\bsW)\in \mathbb R^n \times \mathbb R^{n\times n}: \bsW = \x \t{\x}\}$. For the single pair $(j,k)$, the nonconvex set $X_{jk} = \{(x_j,x_k,w_{jk})\in \mathbb R^3: w_{jk} = x_j x_k\}$ admits a bounded polyhedral relaxation if both variables $x_j,x_k$ are bounded. Assume $x_j\in [\underline x_j, \overline x_j]$ and $x_k \in [\underline x_k, \overline x_k]$. Then, the set $\textrm{conv}(X_{jk}) = X_{jk}^M := \{(x_j,x_k,w_{jk})\in \mathbb R^3\colon
w_{jk} \ge \underline x_k x_j + \underline x_j x_k - \underline x_j \underline x_k, \ 
w_{jk} \ge \overline x_k x_j  + \overline x_j x_k  - \overline x_j \overline x_k, \ 
w_{jk} \le \underline x_k x_j + \overline x_j x_k  - \overline x_j \underline x_k, \ 
w_{jk} \le \overline x_k x_j  + \underline x_j x_k - \underline x_j \overline x_k\}$
is described by the so-called McCormick inequalities \shortcite{mccormick76} and it is proven to provide the tightest relaxation for the set $X_{jk}$ \shortcite{al1983jointly}. Note that the McCormick inequalities specialize to the Fortet reformulation for the product of two binary variables \shortcite{fortet1960applications}. They are also an immediate result of applying the Reformulation-Linearization Technique \shortcite{sherali2013reformulation} to pairs of constraints; for example, the third McCormick inequality is obtained by combining the bounding constraints $x_j - \overline x_j \le 0$ and $x_k - \underline x_k \ge 0$. When independent bounds $[\underline w_{jk}, \overline w_{jk}]$ on $w_{jk}$ are provided, the convex envelope of $X_{jk}$ is, in general, no longer polyhedral \shortcite{anstreicher2021convex,nguyen2018deriving,belotti2010valid}.

However, McCormick inequalities do not yield a tight relaxation of the overall set $X$ as a stronger relationship between the $w_{jk}$ variables is needed. A much tighter relaxation can be obtained by compounding a McCormick-based relaxation with a semidefinite relaxation, \emph{i.e.}, relaxing the matricial constraint $\bsW - \x \t{\x} = 0$ to $\bsW - \x \t{\x} \succeq 0$, which, due to Schur's complement, becomes the PSD constraint
\[
\hat{\bsW}=
\left[
\begin{array}{ll}
1 & \t{\x}\\
\x & \bsW
\end{array}
\right] \succeq 0.
\]
By using either McCormick inequalities or the PSD constraint above, one can obtain two relaxations of [MIQCQP]:
\[
  \begin{array}{llllll}
   \mathbf{(R)} &   \min_{\x,\bsy}       & \textrm{tr}(\bsQ \bsW) + \t{\bsd}\x + \t{\bsr}\bsy\\
                  &\textrm{s.t.} & \textrm{tr}(\bsA_i \bsW) + \t{\bsa_i}\x  + \t{\bsb}\bsy\le c_i & \forall i \in \mathcal{I}\\
                  && \t{\bsf_i}\x + \t{\bsg_i}\bsy \leq h_i & \forall i\ \in \mathcal{I}'\\              &              & (\x,\bsW)\in X_R\\
                  &              & \x\in \mathbb R^n, \bsy \in \mathbb Z^{m},
  \end{array}
\]
by specifying $X_R = \{(\x,\bsW): (x_j,x_k,w_{jk})\in X^M_{jk},\, \forall (j,k)\in N^2\}$ for the McCormick relaxation and $X_R = \big\{(\x,\bsW): \hat{\bsW} \succeq 0 \big\}$ for the SDP relaxation, respectively. Not surprisingly, it was shown that enforcing the SDP relaxation on a McCormick relaxation obtains tighter lower bounds \shortcite{anstreicher2009semidefinite}. This comes with the tradeoff that SDP solvers are not as scalable as Linear Programming (LP) solvers. A workaround has been proposed in \shortciteA{sherali2002enhancing} in the form of {\em SDP cuts}, \emph{i.e.}, linear inequalities that are tangent to the SDP cone: by definition of positive semidefiniteness, for any $\bsa\in \mathbb R^{n+1}$ the linear inequality $\t{\bsa}\hat{\bsW}\bsa\ge 0$ is valid. Such cuts have been adopted both in open-source \shortcite{qualizza2011linear, GonzalezRodriguez2023} and commercial MIQCQP solvers, while other approaches exist that use the SDP relaxation (either solved by an SDP solver or via SDP cuts) within a branch-and-bound approach \shortcite{burer2008finite,buchheim2013semidefinite,dong2016relaxing,gally2018framework,raposaconic}.

Lower bounds for nonconvex MIQCQPs can also be obtained without adding auxiliary variables if all variables are bounded, by adding a PSD quadratic term $\alpha \sum_{j\in N}(x_j - \underline x_j)(\overline x_j - x_j)$, to every quadratic form $\t{\x} \bsQ \x + \t{\bsd}\x$ with $\alpha$ equal to the opposite of the minimum (negative) eigenvalue of~$\bsQ$ \shortcite{adjiman1998global}. Other approaches use reformulations of MIQCQP problems as linear complementarity problems \shortcite{xia2020globally,gondzio2021global}, and special cases such as {\em box-QP}, where variables are only constrained by lower and upper bounds, have been solved through combinatorial optimization methods \shortcite{bonami2018globally}. Due to nonconvex quadratic constraints, solving MIQCQPs requires the use of disjunctions on variables (or combinations thereof) to partition the feasible region. Disjunctions are then used in spatial branch-and-bound \shortcite{elloumi2019global,misener2013glomiqo} as well as lift-and-project and disjunctive programming \shortcite{saxena2010convex,saxena2011convex,dong2018compact}. In general, one can see that the approaches to solve nonconvex MIQCQPs tackle the two types of nonconvexity (integrality in a subset of variables and nonconvex quadratic objective or constraints) by embedding tight lower bounding mechanisms into a branch-and-bound algorithm.

For reasons to be detailed below, in the remainder of this article we focus on the {\em mixed-binary} case where all $\bsy$ variables are binary rather than integer. Apart from streamlining our approach to solving such problems, this comes at little loss of generality, especially given that in most practical applications discrete decisions are indeed represented by binary variables.

\section{Methodology}\label{sec:methodology} 
As highlighted in the previous section, one of the most successful paradigms for solving mixed-integer programs is branch-and-bound. We consider a branch-and-bound framework in which branching decisions are applied to the binary variables $\bsy$, while the continuous variables $\x$ are handled within each node through convex relaxations. In particular, fixing a subset of the binary variables $\bsy$ at a node of the search tree reduces the MIQCQP to a continuous QCQP in the variables $\x$. This resulting problem can then be tackled efficiently using conic relaxations, such as semidefinite or second-order cone relaxations, which provide valid lower bounds. Moreover, when the problem exhibits sparsity, these relaxations can be further decomposed using chordal techniques. This substantially reduces the size of the semidefinite constraints and improves the tractability of the node relaxations within the branch-and-bound procedure.
This strategy is particularly well-suited to problems where the binary decisions determine structural or combinatorial configurations, and the associated continuous variables represent operational or flow decisions. 
The approach enables the exploitation of problem decomposability and structure, especially in applications such as energy networks, transportation design, and facility planning, where the binary decisions encode dispatching or routing choices and the continuous variables determine optimal operations conditioned on those decisions.

The success and the computational efficiency of the branch-and-bound procedure strongly depend on the quality of the relaxation, the early generation of feasible solutions to get good upper bounds, and the branching rules used to define the ensuing subproblems. Exploiting the structure of the relaxation can help speed up the branch-and-bound process by reducing the computational cost at each node. The following sections describe the components of the branch-and-bound algorithm proposed to solve mixed-integer quadratically  constrained quadratic problems formulated as [MIQCQP]. 
Note that, in general, solving nonconvex [MIQCQP] requires spatial branching, \emph{i.e.}, branching on continuous variables. However, we have chosen an approach featuring the following components:
\begin{itemize}
  \item a lower-bounding algorithm based on a convex (conic) continuous relaxation of [MIQCQP];
  \item an efficient, sparsity-based algorithm to obtain lower-bounds from conic relaxations of [MIQCQP];
  \item a heuristic solver that seeks feasible (if suboptimal) MINLP solutions;
  \item a simple branch-and-bound approach where binary variables are branched on.
\end{itemize}
The overarching goal of this article is to show that integrating the above features within a basic branch-and-bound framework can lead to an effective approach for exploring the feasible set of [MIQCQP]. Although combining these components poses nontrivial algorithmic challenges, the resulting method is capable of producing solutions that, while not guaranteed to be optimal, exhibit small optimality gaps.

\subsection{Lower bounds through conic relaxations}
To solve [MIQCQP] efficiently, one would require strong relaxations, particularly in large-scale settings. One common approach to solve such problems is to use MISOCP relaxations, which provide convex approximations of the nonconvex quadratic constraints. However, while MISOCP relaxations can offer good solution bounds, they become computationally prohibitive as the problem size increases due to the large number of conic constraints and their complexity \shortcite{benson2013mixed}. On the other hand, mixed-integer semidefinite programming (MISDP) though also challenging in general, can be highly effective for MIQCQPs that arise from network-structured problems, where they can lead to a more favorable trade-off between solution quality and computational time \shortcite{buchheim2013semidefinite}. Network problems typically exhibit strong sparsity and structural properties that can be exploited using correlative and term sparsity. As we discuss in Section~\ref{sec:sparsity}, we leverage this sparsity structure to focus on MISDP relaxations, since they provide a good balance between solution quality and scalability for large-scale optimization problems. 
Next, we present the SDP relaxation with binary variables of [MIQCQP], which builds upon matrix $\bsW = \x\t{\x}$:

\begin{align*}
\textbf{[MISDP]} \quad \min_{\x, \bsW, \bsy} \quad & \trace(\bsQ \bsW) + \t{\bsq} \x+ \t{\bsr} \bsy \\
\text{s.t.} \quad & \trace(\bsA_i \bsW) + \t{\bsa_i} \x + \t{\bsb_i} \bsy \leq c_i, &\forall i \in \mathcal{I} \\
& \t{\bsf_i}\x + \t{\bsg_i}\bsy \leq h_i, &\forall i\in \mathcal{I'} \\
& \hat{\bsW}= \begin{bmatrix}
1 & \t{\x} \\
{\x} & \bsW
\end{bmatrix}  \\
& \hat{\bsW} \succeq 0, \quad \bsy \in \{0,1\}^m.
\end{align*}


The continuous SDP relaxation of [MISDP], denoted by [SDP-R], is obtained by replacing the binary constraints $\bsy \in \{0,1\}^m$ with $\bsy \in [0,1]^m$. In a branch-and-bound framework, at each node of the tree one solves a node-specific relaxation of [MISDP], where a subset of the $\bsy$ variables is fixed to 0 or 1, while the remaining variables are relaxed to lie in the interval $[0,1]$. Since [MISDP] is a minimization problem, any feasible solution to [MISDP] provides an upper bound on its optimal value, while its continuous relaxation [SDP-R] provides a lower bound. Moreover, because [MISDP] is itself a relaxation of [MIQCQP], its optimal value yields a lower bound for [MIQCQP]. Problem [SDP-R] is convex and can be solved using semidefinite programming techniques. However, the size of the associated SDP matrix ($\hat{\bsW}$), of dimension $(n+1)\times(n+1)$, can become computationally prohibitive as $n$ increases.


\subsection{Exploiting structure: Sparsity}\label{sec:sparsity}
Many real-world problems such as those that arise in energy, gas, and water networks are represented by sparse graphs because the degree of most nodes in these networks is small. We consider two strategies to leverage sparsity based on the definitions of correlative and term sparsity, as described in~\shortciteA{wang2020,Wang:2022}, and implement the general approaches from those papers for [SDP-R]. This allows us to take advantage of any correlative and term sparsity patterns that might be present. Next, we provide more details about the above two types of sparsity. 

There has been a lot of research in the literature on how to reduce the size of SDP constraints, for example by decomposing the matrix $\hat{\bsW}$ into smaller submatrices $\bsW_k$ in a way that 
$\hat{\bsW} \succeq 0$ if and only if all $\bsW_k$ are PSD.
In~\shortciteA{Ghaddar2016, Ghaddar2017}, the so-called \emph{correlative sparsity} (subsets of variables that appear separately in constraints) was considered. Exploiting it allowed solving large-scale water and energy networks. 
Another form of sparsity is \emph{term sparsity}, which is particularly useful when there are few monomials present in the problem. Term sparsity has not yet been applied extensively and, to the best of our knowledge, the only network problem to which it has been applied is the optimal power flow problem~\shortcite{Wang:2022}. Thus, investigating the term sparsity for [SDP-R] is a relevant question as different networks often exhibit a similar structure. The authors in \shortciteA{Wang:2022} combine the two approaches to decompose the SDP matrix of the optimal power flow problem, exploiting the correlative sparsity \shortcite{waki2006sums} and the term sparsity \shortcite{wang2020} at the same time.

Correlative sparsity (CS), introduced in \shortciteA{lasserre2006convergent,waki2006sums,waki2008algorithm}, exploits the fact that each constraint typically involves only a small subset of the decision variables. This structure can be represented using a \emph{correlative sparsity graph}. The vertices of this graph correspond to the continuous variables of the optimization problem, and an edge connects two variables if they appear together in a monomial in the objective or in one of the constraints. This graph is equivalent to the \emph{variable interaction graph} used in sparse SDP formulations. Sparse semidefinite relaxations can be constructed from this graph by decomposing the PSD matrix into smaller matrices associated with subsets of variables. However, to guarantee that the resulting sparse SDP relaxation has the same optimal value as the dense relaxation, the sparsity graph must be chordal. This requirement follows from classical results on positive semidefinite matrix completion \cite{grone1984}, which show that a partial symmetric matrix admits a PSD completion if and only if all its principal submatrices corresponding to the maximal cliques of a chordal graph are positive semidefinite. Therefore, when the correlative sparsity graph is not chordal, a chordal extension is first computed. The goal is to add a set of edges (called fill-in edges) so that the resulting graph becomes chordal, \emph{i.e.}, every cycle of length greater than three contains a chord. Finding the chordal extension with the minimum number of added edges (the minimum fill-in problem) is NP-hard, so practical algorithms rely on efficient heuristics based on vertex elimination orderings. The most widely used is based on the minimum degree heuristic \cite{waki2008algorithm}, which iteratively eliminates the vertex with the smallest degree in the current graph. When a vertex is eliminated, all of its neighboring vertices are connected to each other (introducing fill-in edges if necessary) so that they form a clique. Repeating this process produces a chordal extension of the original graph. The cliques formed during the elimination process correspond to the maximal cliques of the resulting chordal graph, and these cliques determine the blocks of the sparse semidefinite relaxation. Specifically, each maximal clique defines a smaller positive semidefinite matrix involving only the variables in that clique. Since different cliques may share variables, additional equality constraints are imposed to ensure that the overlapping entries of the corresponding matrices agree. This construction yields a sparse SDP relaxation that is equivalent to the dense relaxation while significantly reducing the size of the semidefinite matrices.


Term sparsity (TS), introduced in \shortciteA{wang2020,Wang:2022}, exploits the fact that polynomial constraints typically contain only a small number of monomials. Unlike correlative sparsity, which is based on interactions between variables, term sparsity operates on the monomial basis used in the SDP relaxation. Let $\mathcal{M}$ denote the set of monomials appearing in the objective and constraints of the polynomial optimization problem. In the standard SDP relaxation, the SDP matrix is indexed by monomials from this basis. Term sparsity constructs the \emph{term sparsity graph} by taking the monomials in $\mathcal{M}$ as vertices. Two monomials $m_i,m_j\in\mathcal{M}$ are connected by an edge if their product $m_i m_j$ appears in the polynomial expressions defining the objective or constraints, or if it appears in the construction of the SDP matrix. The connected components of this graph identify subsets of monomials that interact with each other in the relaxation. Instead of constructing a single large SDP matrix indexed by all monomials in $\mathcal{M}$, the SDP relaxation can be decomposed into smaller matrices indexed only by the monomials belonging to each connected component. Since monomials from different components do not interact in the polynomial expressions, the corresponding entries in the SDP matrix are not required. As a result, the SDP relaxation can be written using several smaller semidefinite matrices rather than one large dense matrix.

In this work, we exploit both types of sparsity through the CS-TSSOS framework proposed by \shortciteA{wang2020}. The method first constructs cliques from the correlative sparsity graph to decompose the problem according to variable interactions. It then applies term sparsity within each clique to further partition the monomial basis, yielding smaller SDP blocks. This combined approach significantly reduces the size of the semidefinite matrices while preserving the strength of the relaxation.

To illustrate correlative sparsity further, consider the following polynomial:
$$f(x)=x_1x_4 + x_2x_5 - x_1x_2 - x_4x_5 + x_0(-x_0 + x_1 + x_2 - x_3 + x_4 + x_5).$$

\begin{figure}[!htbp]
    \centering
    \includegraphics[width=0.3\textwidth]{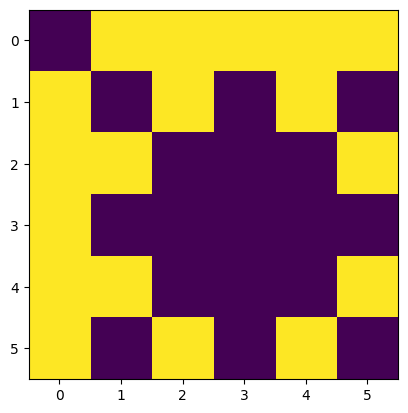}
    \caption{Adjacency matrix of the graph associated with the polynomial example. Light entries indicate variable pairs that appear together in a monomial, while dark entries indicate no direct interaction.}
    \label{fig:example}
\end{figure}

Figure \ref{fig:example} shows the adjacency matrix of the graph representing the variable interactions. Figure \ref{fig:chordalexample}(\subref{subfig:nonchordal}) shows the problem graph where nodes represent variables, $(x_0,\cdots{},x_5)$, and edges indicate that two variables appear together in a term (objective or constraint). This reflects the original sparsity structure of the problem. The graph of the problem is constructed as follows:
\begin{itemize}
    \item $x_0$ and $x_3$ are connected because the monomial $x_0x_3$ appear together in $f$.
    \item  $x_1$ and $x_3$, for example, are not connected as the monomial $x_1x_3$ does not appear in $f$.
\end{itemize}

The graph in Figure~\ref{fig:chordalexample}(\subref{subfig:nonchordal}) is not chordal because it contains the cycle $x_1 - x_2 - x_5 - x_4 - x_1$ of length four without a chord. To apply sparse semidefinite relaxations while preserving equivalence with the dense SDP relaxation, the graph must be chordal. Therefore, we compute a chordal extension by adding edges until every cycle of length greater than three has a chord. In this example, adding the edge $(x_2,x_4)$ triangulates the cycle and yields the chordal graph shown in Figure~\ref{fig:chordalexample}(\subref{subfig:chordal}). This chordal extension is crucial for applying semidefinite programming relaxations efficiently using sparse matrix techniques (\emph{e.g.}, via chordal decomposition). 

\begin{figure}[!htbp]
\centering

\begin{subfigure}{0.45\textwidth}
\centering
\begin{tikzpicture}[scale=1.0,
  every node/.style={circle,draw,thick,minimum size=7mm,font=\small},
  edge/.style={thick},
  chord/.style={thick,dashed}
]
\node (n0) at (0,0) {0};
\foreach \i/\ang in {1/90,2/18,3/-54,4/-126,5/162} {
  \node (n\i) at (\ang:2.4) {\i};
}

\foreach \i in {1,2,3,4,5} \draw[edge] (n0)--(n\i);

\draw[edge] (n1)--(n2);
\draw[edge] (n2)--(n5);
\draw[edge] (n5)--(n4);
\draw[edge] (n4)--(n1);
\end{tikzpicture}

\medskip
\caption{Original non-chordal graph containing the cycle $(x_1,x_2,x_5,x_4)$.}
\label{subfig:nonchordal}
\end{subfigure}
\hfill
\begin{subfigure}{0.45\textwidth}
\centering
\begin{tikzpicture}[scale=1.0,
  every node/.style={circle,draw,thick,minimum size=7mm,font=\small},
  edge/.style={thick},
  chord/.style={thick,dashed}
]
\node (n0) at (0,0) {0};
\foreach \i/\ang in {1/90,2/18,3/-54,4/-126,5/162} {
  \node (n\i) at (\ang:2.4) {\i};
}

\foreach \i in {1,2,3,4,5} \draw[edge] (n0)--(n\i);

\draw[edge] (n1)--(n2);
\draw[edge] (n2)--(n5);
\draw[edge] (n5)--(n4);
\draw[edge] (n4)--(n1);

\draw[chord] (n2)--(n4);
\end{tikzpicture}

\medskip
\caption{Chordal extension obtained by adding the edge $(x_2,x_4)$ (dashed).}
\label{subfig:chordal}
\end{subfigure}

\caption{Polynomial optimization associated graph (a) and chordal extension graph (b). Nodes correspond to variables and edges indicate that two variables appear together in a monomial.}\label{fig:chordalexample}
\end{figure}

The maximal cliques of the chordal graph are
\[
\{x_0,x_3\}, \quad
\{x_0,x_1,x_2,x_4\}, \quad
\{x_0,x_2,x_4,x_5\}.
\]
Using these cliques, the dense SDP matrix can be replaced by three smaller PSD matrices corresponding to these variable subsets. Since the cliques overlap, consistency constraints must be added to ensure that shared entries agree across the matrices. For example, the variables $(x_0,x_2,x_4)$ appear in both of the larger cliques, therefore the corresponding submatrices must coincide on these entries. This decomposition replaces a single dense $6\times6$ SDP matrix with smaller matrices of size at most $4\times4$, reducing computational complexity while preserving the tightness of the relaxation.



\subsection{Upper bounds from a local solver}

In a branch-and-bound algorithm, it is of particular importance to find good feasible solutions as early as possible to increase the fathoming of nodes, reducing the search space. Unfortunately, obtaining feasible solutions for general MIQCQP problems is not an easy task, since it is often already challenging to find feasible solutions for general binary linear programs. In our solution framework for MIQCQP problems, we integrate a MINLP local solver to generate feasible solutions during the branch-and-bound process. 
Several MINLP local solvers can handle integrality as well as nonlinear constraints, and produce feasible solutions or local optima. Global optimality is obviously not required for this purpose, so their use is akin to what primal heuristics do in general mixed-integer optimization solvers.
Given that calls to such solvers are also computationally demanding, the frequency and strategy for invoking this local solver can significantly affect both solution quality and overall computational performance. In general, since we want to find good feasible solutions early on, the local solver is called more often in the initial iterations. The trade-off is clear here: calling the local solver too frequently adds computational overhead, and calling it too infrequently may result in delayed discovery of good, feasible solutions. To balance this trade-off, the local solver is initiated whenever the number of iterations, \texttt{iter}, satisfies the following condition:
\[
\texttt{iter} = 1 \quad \text{or} \quad \texttt{iter} \geq \texttt{run\_local}^{(n_{\texttt{local}} + 3)},
\]
where ${n_\texttt{local}}$ is the number of times the local solver has been previously called and \texttt{run\_local} is a fixed parameter that controls the frequency of the calls: the larger \texttt{run\_local} is, the less frequently the local solver is invoked.

In this framework, when the local solver is called, we provide it with an initial solution derived from the optimal solution of the relaxation at the last solved node. This initial point includes a combination of fixed binary values, enforced by prior branching decisions, and initial values for other binary and continuous variables. Specifically, the binary variables that have not yet been fixed are rounded to the nearest integer, while those constrained by the branching rules are explicitly fixed within the solver. Calling the local solver with different initial solutions at different nodes enables the discovery of different feasible solutions, increasing the likelihood of finding better upper bounds for [MIQCQP] early in the search. 

Apart from the above calls, the local solver is also called at some other points of the algorithm. Recall that the underlying SDP relaxation, [SDP-R], is obtained by relaxing the integrality constraints of the [MISDP] formulation. [SDP-R] is used within the branch-and-bound algorithm to compute a lower bound for [MISDP] and, by extension, for [MIQCQP]. As the branch-and-bound algorithm progresses, whenever an optimal solution of [SDP-R] is also feasible for [MISDP] and improves its upper bound, the local solver is invoked to attempt to obtain a feasible solution for [MIQCQP] close to this point. After these calls, the counter  ${n_\texttt{local}}$ is not updated. The solution of [SDP-R] is passed to the local solver as a starting point but no variable is fixed and, in particular, binary variables $\bsy$ may be modified during the local search. If the local solver finds a feasible solution to [MIQCQP], we subsequently solve [SDP-R] again with all binary variables fixed to the values obtained in this solution. This produces a feasible solution to [MISDP] whose objective value is equal to or better than the one returned by the local solver. As a result, the upper bound for [MISDP] is always at least as strong as the upper bound obtained for [MIQCQP].



\subsection{An SDP-based branch-and-bound algorithm}\label{sec:algorithm}
While no commercial solver can solve mixed-integer optimization problems with SDP constraints, the SCIP-SDP project\footnote{See \url{https://www.opt.tu-darmstadt.de/scipsdp}} uses the open-source SCIP solver to solve MISDP problems \cite{gally2018framework}. It can use several auxiliary SDP solvers as well as SDP cuts, \emph{i.e.}, linear cuts that approximate the SDP cone  \cite{sherali2002enhancing}. Instead of building on top of SCIP-SDP, we have developed our own simple branch-and-bound approach for two reasons: first, to assess our chosen components as a means of finding tight lower bounds, unaffected by other MISDP features; second, to maintain total freedom and simplicity in the algorithmic development, especially the branch-and-bound scheme.
We have developed a solution framework that is generic for [MIQCQP] and in particular for network problems, where sparsity can be exploited. The overall scheme can be seen in Figure~\ref{fig:opfframework}.

\begin{figure}[!htbp]
\begin{tcolorbox}
\begin{tikzpicture}
    \tikzstyle{rect} = [rectangle, draw, text width=3.9cm, text centered, minimum height=1cm]
    \tikzstyle{text-node} = [text width=2cm]
    \tikzstyle{arrow} = [->]
    
    \node[rect] (top) {\textbf{Original \mbox{problem}} \textbf{[MIQCQP]}};
    \node[rect, below=of top] (middle) {\textbf{Mixed-integer SDP relaxation [MISDP]}};
    \node[rect, below=of middle] (bottom) {\textbf{Continuous SDP relaxations [SDP-R]}};
    
    \node[text-node, right = of top, xshift=-0.4cm] (top-right) {Local solver};
    \node[text-node, right = of middle, xshift=-0.4cm] (middle-right) {};
    \node[text-node, right = of bottom, xshift=-0.4cm] (bottom-right) {Conic solver};
    
    \node[right = of top-right, xshift=-0.8cm] (top-right-right) {$UB$ for [MIQCQP]};
    \node[right = of middle-right, xshift=-0.8cm] (middle-right-right) {$UB$ for [MISDP]};
    \node[right = of bottom-right, xshift=-0.8cm] (bottom-right-right) {$LB$ for [MISDP] (and [MIQCQP])};
    
    \draw[arrow] (top) -- (middle);
    \draw[arrow] (middle) -- (bottom);
    
    \draw[decorate,decoration={brace,amplitude=10pt}] 
        ([yshift=-5pt]bottom.south west) -- ([yshift=5pt]middle.north west) 
        node[midway,xshift=-0.8cm,rotate=90] {Branch-and-bound};
\end{tikzpicture}
\end{tcolorbox}
\caption{Diagram of the framework to solve {[MIQCQP]}.}
\label{fig:opfframework}
\end{figure}

We have already introduced all the core ingredients of our branch-and-bound algorithm, whose pseudo-code can be seen in Algorithm~\ref{algorithm:full}. The algorithm terminates when the gap between the current upper and lower bounds of [MISDP] becomes smaller than a predefined tolerance, \texttt{tol}. In this way, the local solver contributes to finding upper bounds, while the branch-and-bound method verifies their quality by tightening the lower bounds through successive refinements of the relaxation. A critical aspect of this framework lies in the quality of the [MISDP] relaxation. Since [MISDP] defines the tightest lower bound achievable within the relaxation, its closeness to the true [MIQCQP] optimum dictates how small the final optimality gap can be. If the [MISDP] solution has significantly weaker lower bound than the [MIQCQP] optimum, the resulting gap may remain large even if high-quality feasible solutions are found. Thus, optimality cannot be guaranteed unless all bounds converge. On the other hand, our framework ensures that we either reach a certified optimal solution of [MIQCQP] (within tolerance) or provide feasible solutions for [MIQCQP] with known optimality gap\textbf{---}the optimality gap defined by the [MISDP] relaxation. This approach is particularly effective for network-structured problems, where sparsity patterns can be exploited to accelerate convergence. We want to emphasize that the gap between the lower bound obtained through the [MISDP] relaxation and the upper bound from a local [MIQCQP] solver does not necessarily converge to zero, even if we fixed all discrete variables $\bsy$: an optimal solution to an SDP relaxation may be even infeasible for [MIQCQP]. The termination criterion for our algorithm is that the [MISDP] relaxation is solved to optimality, and we then compute the [MIQCQP] gap using the best [MIQCQP] solution found.

In the next section, we demonstrate the applicability of the proposed framework on a large-scale real-world problem.


\begin{algorithm}
\caption{Branch-and-Bound Algorithm for [MIQCQP]}
\begin{algorithmic}[1]
\State Set parameters \texttt{tol}, \texttt{timelimit}, \texttt{run\_local}
\State Initialize upper bounds $UB^{\text{MIQCQP}}=UB^{\text{MISDP}} = +\infty$, lower bound $LB = -\infty$
\State Initialize $best^{\text{MIQCQP}}=best^{\text{MISDP}}=\emptyset$
\State Initialize counters: iterations $\texttt{iter} \leftarrow 0$, local solver calls $n_{\texttt{local}} \leftarrow 0$, subproblems $\tau \leftarrow 1$
\State Initialize queue $Q=\{1\}$, $SP^1=$[SDP-R], $LB^1=-\infty$,
\While{\texttt{runtime} < \texttt{timelimit} or $(UB^{\text{MISDP}} - LB)/(\lvert LB \rvert + 10^{-9}) > \texttt{tol}$}
    \State $\texttt{iter} \leftarrow \texttt{iter} + 1$
    \State Choose $k\in Q$ such that $LB^k=\min_{s\in Q} LB^s$. Let $Q \leftarrow Q \backslash \{ k\}$.
    \If{SDP relaxation $SP^k$ has feasible solutions}
        \State Let $(\x^k, \bsy^k)$ be an optimal solution of $SP^k$ and let $obj^k$ be its value.
        \If{$obj^k < UB^{\text{MISDP}}$}
            \If{$(\x^k, \bsy^k)$ is feasible to [MISDP]}
                \State $UB^{\text{MISDP}}\leftarrow obj^k$, $best^{\text{MISDP}}\leftarrow (\x^k, \bsy^k)$
                \State Remove all $s\in Q$ such that $LB^s \geq UB^{\text{MISDP}}$
                \State Call local solver on [MIQCQP] with the solution from $SP^k$ as initial solution
                \State (no binary variable is fixed in this call)
                \If{local solver finds feasible solution better than $UB^{\text{MIQCQP}}$}
                    \State Update $best^{\text{MIQCQP}}$ and $UB^{\text{MIQCQP}}$
                \EndIf
            \Else 
                \State Branch on a binary variable that is fractional in $\bsy^k$, defining $SP^{\tau+1}$, $SP^{\tau+2}$
                \State $Q\leftarrow Q\cup \{ \tau+1, \tau+2\}$, $LB^{\tau+1} \leftarrow obj^k$, $LB^{\tau+2} \leftarrow obj^k$, $\tau\leftarrow\tau+2$
            \EndIf
        \EndIf
    \EndIf
    \If{$\texttt{iter} \geq \texttt{run\_local}^{(n_{\texttt{local}} + 3)}$ or $\texttt{iter} == 1$}
        \State Call local solver on [MIQCQP] with the solution from [SDP-R] as initial solution \break \null \hspace{1cm} (rounding non-fixed binary variables and fixing those already fixed by branching)
        \State $n_{\texttt{local}} \leftarrow n_{\texttt{local}} + 1$
        \If{local solver finds feasible solution better than $UB^{\text{MIQCQP}}$}
            \State Update $best^{\text{MIQCQP}}$ and $UB^{\text{MIQCQP}}$
            \State Solve [SDP-R] with binary variables fixed at the values of the local solver solution
            \If{solution better than $UB^{\text{MISDP}}$ is found}
                \State Update $best^{\text{MISDP}}$ and $UB^{\text{MISDP}}$
            \EndIf
        \EndIf
    \EndIf
    \State $LB\leftarrow \min\{\min_{s\in Q} LB^s, UB^{\text{MISDP}}\}$
\EndWhile
\State \Return $best^{\text{MISDP}}$, $UB^{\text{MISDP}}$, $best^{\text{MIQCQP}}$, $UB^{\text{MIQCQP}}$, $LB$
\end{algorithmic}\label{algorithm:full}
\end{algorithm}

\section{Unit Commitment with Optimal Power Flow Constraints}\label{sec:ucacopf}

The integration of renewable energy sources and increasing demand for resilient and efficient power systems have intensified the need for advanced optimization techniques capable of modeling the inherent complexity of modern electricity markets. In this context, the Unit Commitment (UC) problem, which involves determining the optimal scheduling of generating units, can be used as a model for the reliability and economic operation of power systems. When augmented with AC Optimal Power Flow (ACOPF) constraints, the UC problem captures both operational and network-level physics, resulting in a large-scale, nonconvex MIQCQP formulation.

The nonconvexity of ACOPF constraints can be tackled by convex relaxations, notably SDP relaxations, which offer tight bounds for AC power flow models. Integrating integer decisions from the UC problem yields a MISDP relaxation. Overall, the Unit Commitment problem under OPF constraints, hereafter UC-OPF, combines two optimization problems and forms a MIQCQP formulation of a difficult real-world, large-scale application. The goal of this section is to analyze the potential of the algorithmic approach described in Section~\ref{sec:methodology} by applying it to the UC-OPF problem. We first provide a brief description of the two baseline optimization problems involved in the formulation of UC-OPF.

\subsection{AC optimal power flow problem}

The AC optimal power flow problem, hereafter OPF, was introduced in \shortciteA{Carpentier:1962aa}. Several formulations exist for this class of problems, but the common objective of all of them is to decide how to optimally operate a power distribution network given the nonconvex quadratic constraints imposed by laws of physics and other limits related to the operations. OPF can be adapted to be used for long-term or short-term planning, but regardless of the formulation employed, the use of the power flow equations is what characterizes this problem and also what makes it hard to solve. It is a nonlinear, nonconvex optimization problem and, as mentioned in \shortciteA{Lehmann:2016aa}, it is NP-hard. There is abundant literature devoted to proposing different approaches to solve this problem or similar problems involving power flow constraints: from simplifications of the model using DC power \shortcite{Frank:2012aa} to specific algorithms designed for this particular problem as well as conic based relaxations \shortcite{Lavaei:2012,kocuk2018,savelli2021,bienstock2025,chowdhury2025}. 

\shortciteA{Lavaei:2012} introduce a formulation using rectangular coordinates. We use its terminology and notation; the reader interested in more OPF-related details can find them in this and many other articles in the OPF literature. The problem is modeled as a network of electrical buses (nodes) connected by branches (arcs or edges). The branches, \emph{i.e.}, the connections between the buses, are the representation of transmission lines, cables, transformers, and other equipment. The objective of the system is to transfer electrical energy from the generation buses (supply) to the load buses (demand).

More formally, consider an undirected graph $(N,L)$, where $N$ is the set of buses which contains a subset $G\subset N$ of generator nodes, and a set $L$ of branches $(i,k)$ with $i,k\in N$. Two matrices $\tilde \bsY$ and $\tilde \bsY^{Sh}$ define, respectively, the admittance of the power network while $\tilde Y^{Sh}_{ij}$ is the value of the shunt element at node $i$ for edge $(i,j)$. Parameters $P_i^L$ and $Q_i^L$ are active and reactive power load at node $i\in N$, while $S_{ij}^{\max}$ limits the modulus of the power transfer on $(i,j)\in L$. Admittances $\tilde \bsY$ and $\tilde \bsY^{Sh}$, as well as voltage and power variables in an OPF system, are defined as complex quantities. Let us denote $\Re{z}$ and $\Im{z}$ as the real and imaginary part of $z\in \mathbb C$,
while $\mtr{e}_i$ is the canonical vector in $\mathbb{R}^{\lvert N \rvert}$. The relevant matrices can be defined as follows:
%
 %

{\small
\begin{equation*}
\begin{array}{ll}
\displaystyle \mtr{Y}^G_i = \mtr{e}_i \t{\mtr{e}_i} \Ybt, & \mtr{Y}^B_{ij} = \left( \Yt^{Sh}_{ij} + \Yt_{ij} \right) \mtr{e}_i \t{\mtr{e}_i} - \left( \Yt_{ij} \right) \mtr{e}_i \t{\mtr{e}_j}, \\[3pt]
 \mtr{E}^B_{ij} = \frac{1}{2} \mqty[\Re{\mtr{Y}^B_{ij} + \t{(\mtr{Y}^B_{ij})}} & \Im{\t{(\mtr{Y}^B_{ij})} - \mtr{Y}^B_{ij}} \\ \Im{\mtr{Y}^B_{ij} - \t{(\mtr{Y}^B_{ij})}} & \Re{\mtr{Y}^B_{ij} + \t{(\mtr{Y}^B_{ij})}}], & \mtr{E}^B_{ij} = \frac{1}{2} \mqty[\Re{\mtr{Y}^B_{ij} + \t{(\mtr{Y}^B_{ij})}} & \Im{\t{(\mtr{Y}^B_{ij})} - \mtr{Y}^B_{ij}} \\ \Im{\mtr{Y}^B_{ij} - \t{(\mtr{Y}^B_{ij})}} & \Re{\mtr{Y}^B_{ij} + \t{(\mtr{Y}^B_{ij})}}], \\[10pt]
 \mtr{\bar{E}}^G_i = -\frac{1}{2} \mqty[\Im{\mtr{Y}^G_i + \t{(\mtr{Y}^G_i)}} & \Re{\mtr{Y}^G_i - \t{(\mtr{Y}^G_i)}} \\ \Re{\t{(\mtr{Y}^G_i)} - \mtr{Y}^G_i} & \Im{\mtr{Y}^G_i + \t{(\mtr{Y}^G_i)}}], & \mtr{\bar{E}}^B_{ij} = -\frac{1}{2} \mqty[\Im{\mtr{Y}^B_{ij} + \t{(\mtr{Y}^B_{ij})}} & \Re{\mtr{Y}^B_{ij} - \t{(\mtr{Y}^B_{ij})}} \\ \Re{\t{(\mtr{Y}^B_{ij})} - \mtr{Y}^B_{ij}} & \Im{\mtr{Y}^B_{ij} + \t{(\mtr{Y}^B_{ij})}}], \\[10pt]
 \multicolumn{2}{c}{\mtr{M}_i = \mqty[\mtr{e}_i\t{\mtr{e}_i} & 0 \\ 0 & \mtr{e}_i\t{\mtr{e}_i}] \text{.}}
\end{array}
\end{equation*}
}

Three main classes of variables are defined in an OPF system: voltages $V_i\in \mathbb C$ are defined for each $i\in N$; {\em active} and {\em reactive} power at generator nodes are defined as real scalars $P_i^G$ and $Q_i^G$, respectively. $\tilde P_{ij}$ and $\tilde Q_{ij}$ are, respectively, active and reactive power flows on transmission line $(i,j)\in L$. Voltage variables are arranged in a vector $\bsX\in \mathbb R^{2|N|}$ such that the first $|N|$ elements are the real part of vector $\bsV$ and the last $|N|$ elements are its imaginary part. The objective function, to be minimized, is the sum of quadratic cost functions $C_i(\cdot)$ on generating nodes only. The OPF problem can hence be written as follows:

\begin{subequations}
\begin{spreadlines}{2pt}
\begin{align}
\textbf{[OPF-QP]} \quad \min_{\bm{P}^G, \bm{Q}^G, \bm{\tilde{P}}, \bm{\tilde{Q}}, \bm{X}} \quad& \mathrlap{\sum_{i \in G} C_i(P_i^G)} \notag\\
\text{s.t.}\quad& P_i^G - P_i^L = \trace(\mtr{E}^G_i \var{\mtr{X}} \var{\t{(\mtr{X})}}) & i \in N \label{opf-1}\\
& Q_i^G - Q_i^L = \trace(\mtr{\bar{E}}^G_i \var{\mtr{X}} \var{\t{(\mtr{X})}}) & i \in N\label{opf-2}\\
& P_i^{G, \text{min}} \leq P_i^G \leq P_i^{G, \text{max}} & i \in G \label{opf-3}\\
& Q_i^{G, \text{min}} \leq Q_i^G \leq Q_i^{G, \text{max}} & i \in G \label{opf-4}\\
& (V^{\text{min}}_i)^2 \leq \trace(\mtr{M}_i \var{\mtr{X}} \var{\t{(\mtr{X})}}) \leq (V^{\text{max}}_i)^2 & i \in N & \label{opf-5}\\
& \var{\tilde{P}_{ij}} = \trace(\mtr{E}^B_{ij} \var{\mtr{X}} \var{\t{(\mtr{X})}}) & (i, j) \in L \label{opf-6}\\
& \var{\tilde{Q}_{ij}} = \trace(\mtr{\bar{E}}^B_{ij} \var{\mtr{X}} \var{\t{(\mtr{X})}}) & (i, j) \in L \label{opf-7}\\
& (\var{\tilde{P}_{ij}})^2 + (\var{\tilde{Q}_{ij}})^2 \leq (S^{\text{max}}_{ij})^2 & (i, j) \in L. \label{opf-8}
\end{align}
\end{spreadlines}
\end{subequations}
%
Constraints~\eqref{opf-1}-\eqref{opf-2} enforce nodal power balance at each bus \( i \in N \), ensuring that the net active and reactive power injections (generation minus load) match the total power flows. Constraints~\eqref{opf-3}-\eqref{opf-4} impose operational limits on the active and reactive power outputs of the generators. Voltage magnitude limits at each bus are enforced by constraint~\eqref{opf-5}. Constraints~\eqref{opf-6}-\eqref{opf-7} define the active and reactive power flows on each transmission line \( (i,j) \in L \). Lastly, constraint~\eqref{opf-8} ensures that the apparent power flow on each line does not exceed its thermal limit.
%
We define the feasible set of [OPF-QP] as $F_{\textsc{opf-qp}}=\{(\bm{P}^G, \bm{Q}^G, \bm{\tilde{P}}, \bm{\tilde{Q}}, \bm{X}) \in \mathbb R^{2|G|+2|L|+2|N|}: \eqref{opf-1}-\eqref{opf-8}\}$.

\subsection{Unit commitment problem}

The unit commitment problem, hereafter UC, is a well-known mixed-integer problem, formulated in \shortciteA{Garver:1962} for the first time  
and shown to be an NP-hard problem \shortcite{Lehmann:2016aa}. 
Several attempts to solve this problem efficiently are described in the survey by \shortciteA{Montero:2022}. A significant effort has been devoted to
develop valid cuts to obtain tighter formulations
\shortcite{Ostrowski:2015,Ostrowski:2012,Marcovecchio:2014,Damc-Kurt:2016,Knueven:2018}. 

The UC problem aims at finding an optimal power generation plan for all generator nodes in $G$ in a power network on a planning horizon $T=\{1,2,\ldots{}, \tau\}$ of $\tau$ time instants, to guarantee that the total produced active power satisfies the power demand at each instant. Five classes of variables are defined for each node $i\in N$ and time slot $t\in T$: active and reactive power at generators, $P^G_{it}$ and $Q^G_{it}$ respectively; binary variable $v_{it}$ indicating whether generator $i\in G$ is on at time $t$; and binary variables $y_{it}$ and $z_{it}$, indicating whether generator $i\in G$ was turned on, or turned off, respectively, at time $t$. The UC model is as follows:
\begin{subequations}
\begin{spreadlines}{2pt}
\begin{align}
\textbf{[UC-MIP]}\quad \min_{\bm{v}, \bm{y}, \bm{z}, \bm{P}^G, \bm{Q}^G} \quad& \mathrlap{\sum_{i \in G} \left(\sum_{t \in T}\left(C^2_{it}\var{(P^G_{it})^2}+C^1_{it}\var{P^G_{it}} + C^0_{it} \var{v_{it}} \right) + \sum_{t \in T\backslash\lbrace 0 \rbrace}\left( C^{SU}_{it} \var{y_{it}}\right)\right)} \notag\\
\text{s.t.}\quad& \var{v_{i,t-1}} - \var{v_{it}} + \var{y_{it}} - \var{z_{it}} = 0 &  i \in G, t \in T\backslash \{ 0 \} \label{eq:uc-1}\\
& \var{P^G_{it}} - \var{P^G_{i,t-1}} \leq R^{SU}_{i} \var{v_{i, t-1}} + P^{SU}_{i} \var{y_{it}} & i \in G, t \in T \backslash \{ 0 \} \label{eq:uc-2}\\
& \var{P^G_{i,t-1}} - \var{P^G_{it}} \leq R^{SD}_i \var{v_{it}} + P^{SD}_i \var{z_{it}} & i \in G, t \in T \backslash \{ 0 \} \label{eq:uc-3}\\
& \textstyle{\sum_{k=t-T^U_i+1, \,k \geq 1}^t \var{y_{ik}} \leq \var{v_{it}}} & i \in G, t \in T \label{eq:uc-4}\\
& \textstyle{\var{v_{it}} + \sum_{k=t-T^D_i+1, \,k \geq 1}^t \var{z_{ik}} \leq 1} & i \in G, t \in T\label{eq:uc-5}\\
& P^{G, \text{min}}_i \var{v_{it}} \leq \var{P^G_{it}} \leq P^{G, \text{max}}_i \var{v_{it}} & i \in G, t \in T \label{eq:uc-6}\\
& Q^{G, \text{min}}_i \var{v_{it}} \leq \var{Q^G_{it}} \leq Q^{G, \text{max}}_i \var{v_{it}} & i \in G, t \in T \label{eq:uc-7}\\
& \textstyle{\sum_{i \in G} \var{P^G_{it}} = D_t} & t \in T \label{eq:uc-demand}\\
& \var{y_{it}} + \var{z_{it}} \leq \var{v_{it}} + \var{v_{i,t-1}} & i \in G, t \in T \backslash \{ 0 \} \label{sec:uc-add1}\\
&\var{y_{it}} + \var{z_{it}} \leq 2- (\var{v_{it}} + \var{v_{i, t-1}}) & i \in G, t \in T \backslash \{ 0 \} \label{sec:uc-add2} \\
& v_{it} \in \lbrace 0, 1 \rbrace,\, y_{it} \in \lbrace 0, 1 \rbrace,\, z_{it} \in \lbrace 0, 1 \rbrace & i \in G, t \in T & \text{.}\label{eq:UCtypes}
\end{align}
\end{spreadlines}
\end{subequations}

Constraint \eqref{eq:uc-1} establishes the relation between variables $\bsv$, $\bsy$, and $\bsz$. Constraints \eqref{eq:uc-2} and \eqref{eq:uc-3} limit the increase and decrease in the power generated in a particular generator from one time period to the next, using the ramp-up and ramp-down
rates, $R^{SU}$ and $R^{SD}$, the power generated after being turned on, $P^{SU}$, and the maximum power that the generator can generate before being shut down, $P^{SD}$. Constraints \eqref{eq:uc-4} and \eqref{eq:uc-5} put a minimum duration of the ``on'' or ``off'' status, according to the parameters $T^U$ and $T^D$. Constraints \eqref{eq:uc-6} and \eqref{eq:uc-7} limit the power generated in each generator and constraint \eqref{eq:uc-demand} ensures that demand is met.
Finally, constraints \eqref{sec:uc-add1} and \eqref{sec:uc-add2} are valid cuts
that remove infeasible solutions from the model and, thus, can help reduce the search space of the branch-and-bound scheme.

We define the feasible set of [UC-MIP] as $F_{\textsc{uc-mip}}=\{(\bm{P}^G, \bm{Q}^G, \bsv, \bsy, \bsz) \in \mathbb R^{2|G|\tau + 3|N|\tau}: \eqref{eq:uc-1}-\eqref{eq:UCtypes}\}$. While reactive power variables $\bsQ^G$ play no role in this model (constraint \eqref{eq:uc-7} could be removed and their value could be established after finding an optimal solution to the problem), they are used in the coupled [UCOPF-QP] problem, hence we keep them in the formulation above.

\subsection{Security-constrained unit commitment problem}

Practical, large-scale instances of UC have been tackled successfully to obtain significant savings in the operation of a particular power network \shortcite{Carlson:2012aa}.  However, solving the UC problem might deliver solutions that cannot be used in practice, as they might compromise the stability of the power network via drastic power changes over the time horizon. 
By coupling OPF constraint with an UC formulation, we retain the scheduling capability of UC and guarantee stability by enforcing
the physical restrictions in the power network.

The resulting UC-OPF optimization problem is a challenging one: apart from combining the complexity of both underlying problems, it becomes a large-scale nonconvex MIQCQP due to the extra dimension of the time horizon. 
There has been substantial effort in developing different approaches to obtain feasible solutions for the UC-OPF problem, see for instance the surveys  by \shortciteA{Padhy:2004aa} and \shortciteA{bhardwaj2012unit}.

Most of the effort has focused on using decomposition methods such as Benders decomposition \shortcite{Fu:2006aa, Sifuentes:2007aa, Nasri:2016aa} or outer approximations and Lagrangian relaxations \shortcite{Murillo-Sanchez:1999aa, Ma:1999aa}. Recently, high-quality solutions to realistic industry-scale UC-OPF problems were presented by \shortciteA{brun2025} as part of the Grid Optimization Competition \shortcite{elbert2024}, who decompose the problem and apply a penalty-based alternating direction method to solve it. Similarly, the work of \shortciteA{zhang2023} employs an alternating direction algorithm to efficiently obtain good feasible solutions for the UC-OPF problem. Yet, in general, most of the proposed solution methods do not provide optimality guarantees for the obtained solutions, not even for local optimality in some cases. 

Building upon the SDP relaxation for the OPF problem, we introduce a new framework to solve the UC-OPF problem. We employ the SDP relaxation, adapting it to the UC-OPF case, and then embed the resulting mixed-integer SDP problem into the branch-and-bound scheme described in Section~\ref{sec:algorithm} to find an optimal solution to the relaxation. This optimal solution is then used to certify the optimality of a solution to the UC-OPF problem obtained by the use of a local solver for MINLP problems or, at least, to provide a valid optimality gap for it.

The structure of the combined UC-OPF problem presented below is straightforward: we combine the variables and constraints from both [OPF-QP] and [UC-MIP]. The main difference here is that, for every time period in the UC problem, we now have an entire OPF problem associated with it, resulting in a challenging multi-period nonlinear problem. Note that constraint \eqref{eq:uc-demand} can be removed, as it is subsumed by the power load constraints \eqref{opf-1} and \eqref{opf-2}. 
The resulting model is as follows:
\begin{subequations}
\label{eq:ucacopf-original}
\begin{align}
\textbf{[UCOPF-QP]}\min_{\bm{v}, \bm{y}, \bm{z}, \bm{P}, \bm{Q}, \bm{\tilde{P}}, \bm{\tilde{Q}}, \bm{X}}\quad& \mathrlap{\sum_{i \in G} \left(\sum_{t \in T}\left(C^2_{it}\var{(P^G_{it})^2}+C^1_{it}\var{P^G_{it}} + C^0_{it} \var{v_{it}} \right) + \sum_{t \in T\backslash\lbrace 0 \rbrace}\left( C^{SU}_{it} \var{y_{it}}\right)\right)} \notag\\
\text{s.t.}\quad
& (\bsv, \bsy, \bsz, \bsP^G, \bsQ^G) \in F_{\textsc{uc-mip}}\\
& (\bsP_t, \bsQ_t, \tilde \bsP_t, \tilde \bsQ_t, \bsX_t) \in F_{\textsc{opf-qp}}
& \forall t \in T,\qquad\qquad\qquad
\end{align}
\end{subequations}
where all variables of [OPF-QP] now have an extra instant index $t$ and vectors $\bsP^G$ and $\bsQ^G$ are subvectors of $\bsP$ and $\bsQ$, respectively.

We obtain a MISDP relaxation by replacing the nonlinear, nonconvex terms $\bsX \t{\bsX}$ by $\bsW$ and by adding the extra constraint $\bsW \succeq 0$. Let us define the feasible set of the SDP relaxation of [OPF-QP] by replacing constraints \eqref{opf-1}, \eqref{opf-2}, 
\eqref{opf-5}, \eqref{opf-6}, and \eqref{opf-7} with constraints

\begin{subequations}
\begin{align}
& P_i^G - P_i^L = \trace(\bsE_i^G \bsW_t) &  i \in N, t\in T\label{eq:ActPowSDP}\\
& Q_i^G - Q_i^L = \trace(\bar \bsE_i^G \bsW_t) &  i \in N, t\in T \label{eq:ReaPowSDP}\\
&(V_i^{\min})^2 \le \trace(\bsM_i \bsW_t) \le (V_i^{\max})^2 &  i \in N, t\in T \label{eq:modRangeSDP}\\
& \tilde P_{ijt} = \trace(\bsE_{ij}^B \bsW_t) & (i,j)\in L, t\in T\label{eq:defActPowTransferSDP}\\
& \tilde Q_{ijt} = \trace(\bar \bsE_{ij}^B \bsW_t) & (i,j)\in L, t\in T \label{eq:defReaPowTransferSDP}\\
& \bsW_t \succeq 0 &  t\in T.\label{eq:Wpsd} 
\end{align}
\end{subequations}
Note that the resulting formulation is not a pure MISDP, as quadratic terms in the objective and constraint \eqref{opf-8} are handled separately via conic reformulations as done by \cite{Lavaei:2012}. We now apply the approach presented in Section~\ref{sec:methodology} to solve [UCOPF-QP]. This is a mixed-integer nonlinear nonconvex optimization problem, as we have binary variables from the unit commitment part, a quadratic objective function, and quadratic constraints related to the power equations.

Let us first define the feasible set of the SDP relaxation: $F_{\textsc{opf-sdp}}=\{(\bm{P}^G, \bm{Q}^G, \bm{\tilde{P}}, \bm{\tilde{Q}}, \bm{W}) \in \mathbb R^{2|G|+2|L|+2|N|}: \eqref{opf-3}, \eqref{opf-4}, \eqref{opf-8}, \eqref{eq:ActPowSDP}, \eqref{eq:ReaPowSDP}, \eqref{eq:modRangeSDP}, \eqref{eq:defActPowTransferSDP}, \eqref{eq:defReaPowTransferSDP}, \eqref{eq:Wpsd} \}$.
Then the [UCOPF-MISDP] relaxation can be written as follows:

\begin{subequations}
\label{eq:ucacopf-sdprelaxation}
\begin{align}
\textbf{[UCOPF-MISDP]} & \notag \\
\min_{\bm{v}, \bm{y}, \bm{z}, \bm{P}, \bm{Q}, \bm{\tilde{P}}, \bm{\tilde{Q}},  \bm{W}}\quad& \mathrlap{\sum_{i \in G} \left(\sum_{t \in T}\left(C^2_{it}\var{(P^G_{it})^2}+C^1_{it}\var{P^G_{it}} + C^0_{it} \var{v_{it}} \right) + \sum_{t \in T\backslash\lbrace 0 \rbrace}\left( C^{SU}_{it} \var{y_{it}}\right)\right)} \notag\\
\text{s.t.}\qquad
& (\bsv, \bsy, \bsz, \bsP^G, \bsQ^G) \in F_{\textsc{uc-mip}}\\
& (\bsP_t, \bsQ_t, \tilde \bsP_t, \tilde \bsQ_t, \bsW_t) \in F_{\textsc{opf-sdp}} & \forall t \in T.
\end{align}
\end{subequations}


Problem [UCOPF-QP] is then a suitable test case for the methodology described in Section~\ref{sec:methodology} and, more precisely, Algorithm~\ref{algorithm:full}. With respect to the scheme of Algorithm~\ref{algorithm:full}, the [MIQCQP] problem and its SDP relaxation [MISDP] now become [UCOPF-QP] and [UCOPF-MISDP], respectively.


\subsection{Instances for the case study}
We construct the instances for our computational study by combining instances for the OPF problem and instances from the UC problem. In order to do so, we take three well-known OPF instances from MATPOWER \shortcite{zimmerman2010matpower}: \texttt{case6ww}, \texttt{case24\_ieee\_rst}\footnote{This instance has multiple generators at its buses, which would require extra decision variables in the OPF model. We do take this into account in our experiments by defining $P^G_{i}$ as the sum of the power of all generators at node $i$; however,  for the sake of simplicity we have chosen to exclude the corresponding extra constraints from the description of our models in this article.}, and \texttt{case118}, which correspond to small and medium-size instances. On top of these OPF instances we add the unit commitment data from \shortciteA{castillo2016}. It is worth noting that there is some overlap between the OPF parameters in MATPOWER and in \shortciteA{castillo2016}, but they are slightly different to one another. Thus, we have to make some modifications to the combined instances to ensure feasibility of the resulting UC-OPF instances. Hereafter, we use Case-6, Case-24 and Case-118 to refer to these combined instances. Table~\ref{tab:cases} contains information regarding the sizes of the underlying power networks. These are the instances used in Section \ref{sec:results}. If sparsity is not exploited, these instances have SDP matrices of size twice the number of buses ($2\cdot\lvert N \rvert$) and, as we discuss in Section~\ref{sec:results}, even the SDP relaxation for the OPF problem alone is not solvable for Case-118 without exploiting the structure.

\begin{table}[!htbp]
    \centering
    \begin{tabular}{l|ccc}
        \hline
        {Case} &  Buses& Generators & Branches \\
        \hline \hline
        Case-6 & 6 & 3 & 11\\
        Case-24& 24  & 32 & 38  \\
        Case-118 & 118  & 54 & 186\\
        \hline
    \end{tabular}
    \caption{Sizes of baseline OPF instances. Branches account for transmission lines and transformers.}
    \label{tab:cases}
\end{table}

\subsection{Exploiting sparsity}
\label{sec:ucacopf-sparsity}


As mentioned in the previous sections, the SDP relaxation for the OPF problem on its own is very expensive to solve even for medium-scale instances and, for UC-OPF, which has SDP constraints for the different time periods, complexity increases even further. To overcome the challenge of solving the SDP relaxation, we exploit the sparsity of the power network and replace the semidefinite constraints with smaller semidefinite constraints based on the correlative and term sparsity presented in Section \ref{sec:sparsity}.  Matrix $\Wb$ contains all possible monomials of degree $2$ that can appear in [UCOPF-SDP]. As defined, the SDP constraint $\Wb \succeq 0$ scales in size with the increase of the number of nodes in the problem. This complicates its use for large problems.

\begin{figure}[!htbp]
    \begin{subfigure}[b]{0.33\textwidth}
        \centering
    \includegraphics[width=\textwidth]{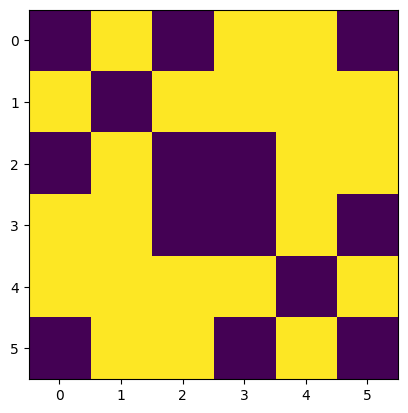}
        \caption{Case-6}
        \vspace*{4mm}
        \label{fig:case6}
    \end{subfigure}
    \begin{subfigure}[b]{0.33\textwidth}
        \centering
         \includegraphics[width=\textwidth]{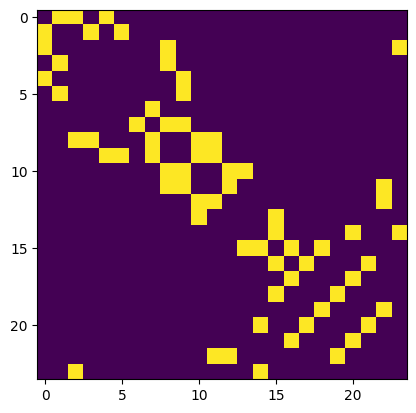}
        \caption{Case-24}
        \vspace*{4mm}
        \label{fig:case24}
    \end{subfigure}    
    \begin{subfigure}[b]{0.33\textwidth}
        \centering
    \includegraphics[width=\textwidth]{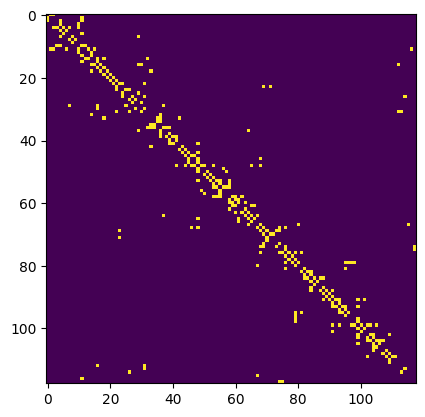}
        \caption{Case-118}
        \vspace*{4mm}
        \label{fig:case118}
    \end{subfigure}
    \caption{Sparsity of the [OPF-QP] problem. Light entries indicate variable pairs that appear together in a monomial, while dark entries indicate no direct
interaction.}
    \label{fig:sparsity}
\end{figure}

In this case, both correlative and term sparsity methods deliver smaller matrices (blocks) that can be used to decompose the SDP constraints.  Note that this needs to be done for a single $\Wb$ and the rest of the $\Wb_t$ will follow the same block structure since the $\Wb$ variables at different time periods are identical in terms of the sparsity pattern. This is done once at the root node of the branch-and-bound tree, as the same pattern follows at each node. Figure \ref{fig:sparsity} shows the network sparsity patterns, which become more pronounced as the network size increases. 

Once the sparsity structure is exploited, the maximum size of the resulting SDP matrices is notably smaller than the size of the original ones. Table 2 reports the maximum SDP block sizes and total number of blocks for the different cases, showing that size reduction becomes more pronounced as the number of buses increases. The number in parentheses in the Max SDP size column indicates how many SDP matrices have that size.

\begin{table}[!htbp]
    \centering
    \begin{tabular}{l|c|cc}
        \hline
        {Case} &  Original SDP size& Max SDP size & \#SDP Blocks \\
        \hline \hline
        Case-6 & 12 & 8 (4) & 10\\
        Case-24& 48  & 11 (1) & 50  \\
        Case-118 &236  &11 (6)& 266\\
        \hline
    \end{tabular}
    \caption{The size of the SDP matrices for the original problem and exploiting sparsity.}
    \label{table:sparsity}
\end{table}

For instance, for Case-118, the SDP matrix is reduced from 236$\times$236 to 266 SDP matrices of smaller dimension with a maximum size of 11$\times$11. As the interior point method of solving the SDP relaxation is significantly impacted by the largest SDP matrix size, this reduces the computational time significantly for such instances, as seen in the computational results.

\begin{figure}[!htbp]
        \centering
    \includegraphics[width=0.6\textwidth]{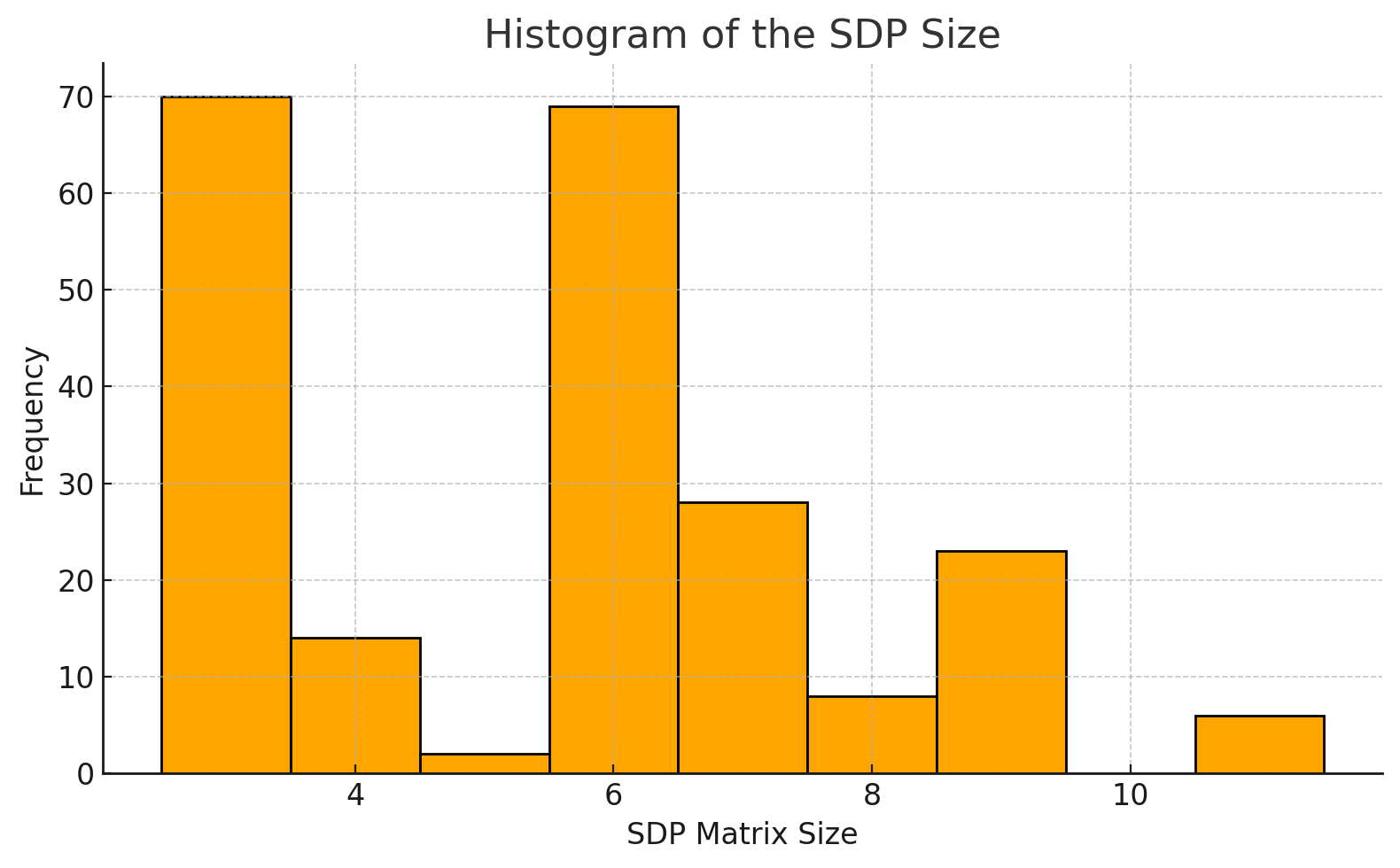}
        \caption{Histogram of the resulting sizes of SDP matrices in Case-118.}
        \vspace*{4mm}
        \label{fig:histogramcase118}
    \end{figure}

Figure \ref{fig:histogramcase118} presents a histogram of the SDP sizes of the matrices resulting from exploiting problem sparsity for Case-118. The frequency distribution of these sizes ranges from 3 to 11. We can clearly observe clusters around 3 and 6 as the most frequent values and secondary modes around 7 and 9.

\section{Numerical Results for the UC-OPF Problem}\label{sec:results}

\subsection{Computational setting}

\subsubsection{Running environment}
All the executions reported in this paper have been run on the supercomputer Finisterrae~III, at the Galicia Supercomputing Centre (CESGA). Specifically, we used computational nodes powered with two thirty-two-core Intel Xeon Ice Lake 8352Y CPUs with 256GB of RAM and 1TB SSD. Algorithm~\ref{algorithm:full} has been implemented in Python~\shortcite{python}.

\subsubsection{Test instances}
In order to enrich the computational analysis, we define a series of variations of the baseline instances, Case-6, Case-24, and Case-118, as follows:
\begin{description}
    \item[Initial conditions.] We duplicate the number of instances by defining variations of them in which the initial conditions of the generators at $t=0$ are set to 0. We denote the original instances by Case-6a, Case-24a, and Case-118a and the new ones by Case-6b, Case-24b, and Case-118b.
    \item[Random perturbations.] For each of the above 6 instances we create another 6 variations in which the demand and bounds on the reactive power at the nodes are randomly perturbed with average noise levels of 2\%, 4\% and 6\% (perturbations are independent across buses and across time periods). Two instances are generated for each of these three levels of noise, resulting in a total of 36 new instances. 
    \item[Time periods.] For each instance we work with three different values for the time horizon: 4 periods, 12 periods, and the full 24 periods. 
\end{description}

This process yields a total of 126 instances, obtained as the product of $6$ baseline instances, $7$ variations, and $3$ time horizons. The test instances are publicly available and can be downloaded at \url{https://bit.ly/data-UC-ACOPF}.

\subsubsection{Algorithm configuration}
Our framework consists of four main components: the local solver, the SDP solver, the sparsity exploitation, and the branch-and-bound algorithm in charge of managing the binary variables. Below we describe in detail how each of these elements has been set up for the numerical analysis\footnote{In our implementation, equations \eqref{eq:ActPowSDP} and \eqref{eq:ReaPowSDP} are included as inequality constraints. Given the monotonically non-decreasing objective function, this yields an equivalent formulation while making feasible solutions easier to find, which may improve algorithmic performance.}.

\begin{description}
    \item[Local solver.] We begin by describing the local solver. We use Knitro \shortcite{Byrd:2006} as the nonlinear programming solver, interfaced through Pyomo \shortcite{Bynum:2021aa, Hart:2011aa}. We rely primarily on Knitro's default settings, with a few modifications to better suit our iterative node-based strategy. Specifically, we set \texttt{mip\char`_multistart} to 1, \texttt{convex} to 0, and \texttt{mip\char`_terminate} to 1. The last setting instructs Knitro to return a feasible solution promptly, without spending additional time trying to refine it. To safeguard against excessive runtimes in difficult nodes, we also cap the solver time limit using the \texttt{mip\char`_maxtime\char`_real} option. 
    \item[SDP solver.]  We use conic programming solver MOSEK \shortcite{mosek}. We formulate the SDP relaxation using MOSEK’s Python-based Fusion API, which allows flexible modeling of conic problems and provides direct access to advanced solver features \shortcite{mosek}. 
    \item[Sparsity exploitation] We follow the approach described in \shortciteA{wang2020,Wang:2022}, to maintain the theoretical properties of the decompositions. We use their Julia \shortcite{Bezanson:2017} package to obtain the decomposition of the $\Wb_t$ matrices. Their approach is called CS-TSSOS and, as discussed in Section~\ref{sec:ucacopf-sparsity}, it combines correlative sparsity \shortcite{waki2006sums} and term sparsity \shortcite{wang2020}.
    \item[Branch-and-bound algorithm.] As mentioned above, we have implemented from scratch a branch-and-bound algorithm in Python, tailored for solving mixed-integer QCQP problems. The parameters of Algorithm~\ref{algorithm:full} are set up as follows:
    \begin{itemize}
        \item \texttt{tol}. Stopping tolerance is set to $0.02$, \emph{i.e.}, $2\%$ optimality gap.
        \item \texttt{timelimit}. Time limit is set to 3600~seconds (1~hour).
        \item \texttt{run\_local}. All instances are solved 4 times, each with a different value of \texttt{run\_local} parameter: $1.1$, $1.25$, $1.5$, and $2$. The larger the value of \texttt{run\_local}, the smaller the number of calls to the local solver.
    \end{itemize}
\end{description}

\subsection{First results (without accounting for sparsity)}
We first present results for [UCOPF-QP] solved without exploiting any sparsity. Table~\ref{tab:first_results} summarizes the outcomes for each baseline instance and number of periods.
\begin{table}[!htbp]
\centering
{\small
\begin{tabular}[t]{lcccrcrc}
\toprule
\multicolumn{1}{c}{Instance} & \multicolumn{1}{c}{Periods} & \multicolumn{1}{c}{MISDP-Gap} & \multicolumn{1}{c}{MIQCQP-Gap} & \multicolumn{1}{c}{Time} & \multicolumn{1}{c}{LS\,\%} & \multicolumn{1}{c}{Iter} & \multicolumn{1}{c}{UB jumps}\\
\midrule
\cellcolor{gray!25}{Case-6a} & \cellcolor{gray!25}{4} & \cellcolor{gray!25}{1.05\,\%} & \cellcolor{gray!25}{1.17\,\%} & \cellcolor{gray!25}{3.45} & \cellcolor{gray!25}{61.72\,\%} & \cellcolor{gray!25}{23.9} & \cellcolor{gray!25}{2.0}\\
Case-6a & 12 & 0.99\,\% & 2.11\,\% & 45.44 & 39.24\,\% & 187.3 & 2.3\\
\cellcolor{gray!25}{Case-6a} & \cellcolor{gray!25}{24} & \cellcolor{gray!25}{1.90\,\%} & \cellcolor{gray!25}{1.90\,\%} & \cellcolor{gray!25}{119.73} & \cellcolor{gray!25}{35.67\,\%} & \cellcolor{gray!25}{338.1} & \cellcolor{gray!25}{2.7}\\
Case-6b & 4 & 1.70\,\% & 3.11\,\% & 13.29 & 67.05\,\% & 59.6 & 2.4\\
\cellcolor{gray!25}{Case-6b} & \cellcolor{gray!25}{12} & \cellcolor{gray!25}{1.75\,\%} & \cellcolor{gray!25}{1.99\,\%} & \cellcolor{gray!25}{117.34} & \cellcolor{gray!25}{15.63\,\%} & \cellcolor{gray!25}{765.6} & \cellcolor{gray!25}{2.9}\\
Case-6b & 24 & 1.98\,\% & 1.98\,\% & 402.11 & 17.94\,\% & 1440.4 & 3.6\\
\cellcolor{gray!25}{Case-24a} & \cellcolor{gray!25}{4} & \cellcolor{gray!25}{1.88\,\%} & \cellcolor{gray!25}{1.89\,\%} & \cellcolor{gray!25}{1994.13} & \cellcolor{gray!25}{11.95\,\%} & \cellcolor{gray!25}{1272.7} & \cellcolor{gray!25}{6.9}\\
Case-24a & 12 & 1.90\,\% & 1.90\,\% & 1572.03 & 43.63\,\% & 163.1 & 4.3\\
\cellcolor{gray!25}{Case-24a} & \cellcolor{gray!25}{24} & \cellcolor{gray!25}{4.88\,\%} & \cellcolor{gray!25}{4.88\,\%} & \cellcolor{gray!25}{3462.28} & \cellcolor{gray!25}{76.17\,\%} & \cellcolor{gray!25}{80.1} & \cellcolor{gray!25}{3.1}\\
Case-24b & 4 & 14.87\,\% & 14.88\,\% & \textbf{---} & 8.15\,\% & 2315.3 & 6.6\\
\cellcolor{gray!25}{Case-24b} & \cellcolor{gray!25}{12} & \cellcolor{gray!25}{5.09\,\%} & \cellcolor{gray!25}{5.09\,\%} & \cellcolor{gray!25}{\textbf{---}} & \cellcolor{gray!25}{31.39\,\%} & \cellcolor{gray!25}{439.9} & \cellcolor{gray!25}{4.0}\\
Case-24b & 24 & 7.94\,\% & 7.94\,\% & \textbf{---} & 61.77\,\% & 151.0 & 3.0\\
\bottomrule
\end{tabular}
}
\caption{\label{tab:first_results}Average results by instance and number of periods (no sparsity exploitation).}
\end{table}

Recall that, for each baseline instance there are a total of 7 variations of it (the original one plus 6 random perturbations) and that each of them was executed with 4 different values of parameter \texttt{run\_local}. The numbers reported in Table~\ref{tab:first_results} are computed as follows. For each of the 7 variations of an instance, we take the best-performing execution among the 4 different values for \texttt{run\_local}, according to the average between the gap for MISDP and the gap for MIQCQP. Once we have the 7 best-performing executions, we take the average of the different values of interest, which are then collected in Table~\ref{tab:first_results}. Column ``LS\,\%'' shows the average percentage of the overall running time spent by the local solver and column ``UB jumps'' contains the average number of times that a new ``best'' primal solution was found by the algorithm. Note that, when presenting the results grouped as in Table~\ref{tab:first_results}, it is not possible to assess the impact of the choices for parameter \texttt{run\_local}. This impact is  discussed in Section~\ref{sec:local}, where we present results disaggregated by the value of \texttt{run\_local}.

Table~\ref{tab:first_results} already shows some patterns that will be common in all ensuing executions. Instances with the initial conditions set to 0 (``b'' instances) tend to be harder to solve. Also, as expected, a larger number of buses and a larger number of periods also make instances harder to solve. Indeed, it is worth noting that, without exploiting sparsity, no Case-118 instance was solved within the time limit and, more importantly, not even a feasible solution of [MISDP] was found, which is the reason why no result is reported for them in Table~\ref{tab:first_results}.

Maybe the most important aspect to note in this first set of results is the tightness of the MISDP relaxations, which was one of the main motivations behind the algorithmic approach we have developed in this paper. In order to assess this claim one just needs to observe that the difference between MISDP-Gap and MIQCQP-Gap is very small in almost all the instances. Recall that, whenever a new ``best'' solution is found for [MISDP], the algorithm calls the local solver to obtain a close solution for [MIQCQP]. The small difference between MISDP-Gap and MIQCQP-Gap indicates that these calls are successful at finding solutions of [MIQCQP] whose objective values are very close to those obtained for the original [MISDP] ones.

\subsection{Main results: Exploiting sparsity}

Table~\ref{tab:main_results} is analogous to Table~\ref{tab:first_results}, but now all the executions have been made accounting for the sparsity structure in MIQCQP as discussed in Section~\ref{sec:ucacopf-sparsity}. The results in both tables are qualitatively very similar, but a closer look at optimality gaps and running times shows that, although the performance is very similar for Case-6 instances, accounting for sparsity already starts to pay off in the Case-24 ones. This is not surprising, since we already showed in Table~\ref{table:sparsity} that the size reduction on the SDP matrices increases with the number of buses.

\begin{table}[!htbp]
\centering
{\small
\begin{tabular}[t]{lcccrcrc}
\toprule
\multicolumn{1}{c}{Instance} & \multicolumn{1}{c}{Periods} & \multicolumn{1}{c}{MISDP-Gap} & \multicolumn{1}{c}{MIQCQP-Gap} & \multicolumn{1}{c}{Time} & \multicolumn{1}{c}{LS\,\%} & \multicolumn{1}{c}{Iter} & \multicolumn{1}{c}{UB jumps}\\
\midrule
\cellcolor{gray!25}{Case-6a} & \cellcolor{gray!25}{4} & \cellcolor{gray!25}{1.05\,\%} & \cellcolor{gray!25}{1.17\,\%} & \cellcolor{gray!25}{3.60} & \cellcolor{gray!25}{55.39\,\%} & \cellcolor{gray!25}{24.9} & \cellcolor{gray!25}{1.7}\\
Case-6a & 12 & 0.99\,\% & 2.07\,\% & 43.30 & 29.42\,\% & 188.3 & 2.4\\
\cellcolor{gray!25}{Case-6a} & \cellcolor{gray!25}{24} & \cellcolor{gray!25}{1.90\,\%} & \cellcolor{gray!25}{1.90\,\%} & \cellcolor{gray!25}{173.74} & \cellcolor{gray!25}{38.16\,\%} & \cellcolor{gray!25}{339.7} & \cellcolor{gray!25}{3.0}\\
Case-6b & 4 & 1.70\,\% & 3.11\,\% & 15.27 & 67.64\,\% & 60.6 & 2.6\\
\cellcolor{gray!25}{Case-6b} & \cellcolor{gray!25}{12} & \cellcolor{gray!25}{1.75\,\%} & \cellcolor{gray!25}{1.99\,\%} & \cellcolor{gray!25}{152.04} & \cellcolor{gray!25}{14.90\,\%} & \cellcolor{gray!25}{774.0} & \cellcolor{gray!25}{2.9}\\
Case-6b & 24 & 1.97\,\% & 1.97\,\% & 449.06 & 8.67\,\% & 1458.0 & 3.0\\
\cellcolor{gray!25}{Case-24a} & \cellcolor{gray!25}{4} & \cellcolor{gray!25}{1.84\,\%} & \cellcolor{gray!25}{1.84\,\%} & \cellcolor{gray!25}{1063.67} & \cellcolor{gray!25}{27.18\,\%} & \cellcolor{gray!25}{1922.4} & \cellcolor{gray!25}{6.3}\\
Case-24a & 12 & 1.88\,\% & 1.88\,\% & 1326.12 & 75.67\,\% & 450.3 & 4.7\\
\cellcolor{gray!25}{Case-24a} & \cellcolor{gray!25}{24} & \cellcolor{gray!25}{3.33\,\%} & \cellcolor{gray!25}{3.33\,\%} & \cellcolor{gray!25}{3309.43} & \cellcolor{gray!25}{90.99\,\%} & \cellcolor{gray!25}{127.7} & \cellcolor{gray!25}{3.4}\\
Case-24b & 4 & 7.39\,\% & 7.41\,\% & \textbf{---} & 8.38\,\% & 6551.0 & 6.4\\
\cellcolor{gray!25}{Case-24b} & \cellcolor{gray!25}{12} & \cellcolor{gray!25}{3.91\,\%} & \cellcolor{gray!25}{3.91\,\%} & \cellcolor{gray!25}{\textbf{---}} & \cellcolor{gray!25}{48.57\,\%} & \cellcolor{gray!25}{1930.3} & \cellcolor{gray!25}{5.3}\\
Case-24b & 24 & 9.45\,\% & 9.45\,\% & \textbf{---} & 64.04\,\% & 686.3 & 2.9\\
\cellcolor{gray!25}{Case-118a} & \cellcolor{gray!25}{4} & \cellcolor{gray!25}{0.01\,\%} & \cellcolor{gray!25}{0.01\,\%} & \cellcolor{gray!25}{17.40} & \cellcolor{gray!25}{85.34\,\%} & \cellcolor{gray!25}{1.0} & \cellcolor{gray!25}{1.0}\\
Case-118a & 12 & 0.06\,\% & 0.06\,\% & 56.84 & 87.03\,\% & 1.0 & 1.0\\
\cellcolor{gray!25}{Case-118a} & \cellcolor{gray!25}{24} & \cellcolor{gray!25}{0.04\,\%} & \cellcolor{gray!25}{0.04\,\%} & \cellcolor{gray!25}{137.84} & \cellcolor{gray!25}{86.51\,\%} & \cellcolor{gray!25}{1.0} & \cellcolor{gray!25}{1.0}\\
Case-118b & 4 & 0.00\,\% & 0.00\,\% & 13.32 & 82.08\,\% & 1.0 & 1.0\\
\cellcolor{gray!25}{Case-118b} & \cellcolor{gray!25}{12} & \cellcolor{gray!25}{0.02\,\%} & \cellcolor{gray!25}{0.02\,\%} & \cellcolor{gray!25}{45.58} & \cellcolor{gray!25}{84.75\,\%} & \cellcolor{gray!25}{1.0} & \cellcolor{gray!25}{1.0}\\
Case-118b & 24 & 0.09\,\% & 0.09\,\% & 126.68 & 85.71\,\% & 1.0 & 1.0\\
\bottomrule
\end{tabular}
}
\caption{\label{tab:main_results}Average results by instance and number of periods (accounting for sparsity).}
\end{table}

Importantly, Table~\ref{tab:main_results} also shows that Case-118 instances can be successfully solved once sparsity is accounted for. Indeed, a closer look at the results shows that they are solved too easily: the optimality gaps are very small and, surprisingly, both the number of iterations and UB jumps take value 1. After a closer inspection at the generator's data taken from MATPOWER \shortcite{zimmerman2010matpower}, we observed that, although Case-6 and Case-24 instances had strictly positive values for $P^{\text{min}}_i$, all these values are set to $0$ for Case-118. This implies that there is no minimum generation needed at the generators which is one of the main complexities behind the combinatorial nature of the UC problem. These values for $P^{\text{min}}_i$, combined with the rest of the data of the problem, result in the instance being solved at the root node.

In order to get a more informative analysis of the Case-118 instances, in Section~\ref{sec:case118finer} we enhance the analysis for these instances by using the $P^{\text{min}}_i$ data from \shortciteA{castillo2016}, along with the definition of additional Case-118 instances by controlled variations of these baseline $P^{\text{min}}_i$ values.

\begin{figure}[!htbp]
\centering
\includegraphics[width=\textwidth]{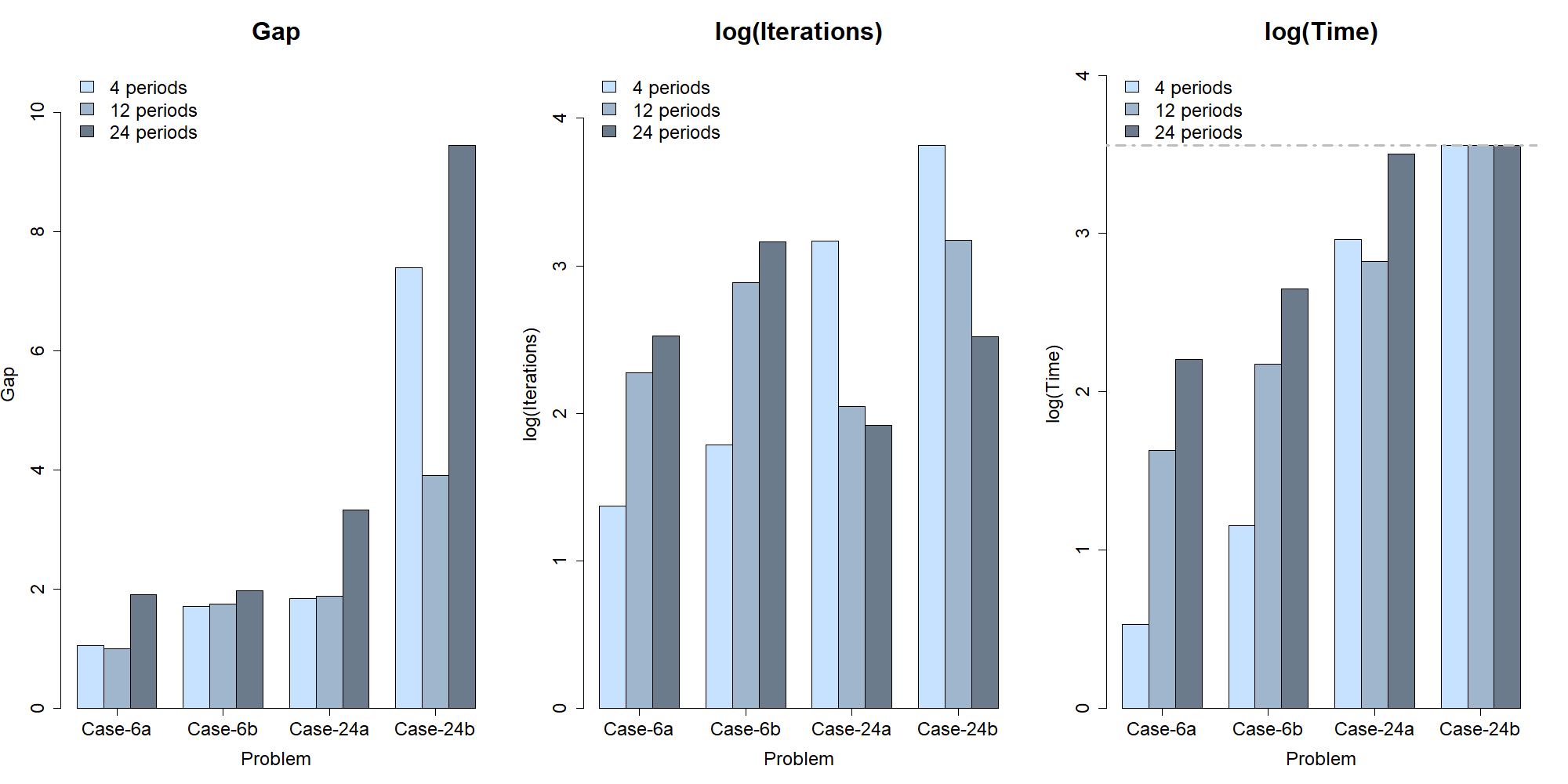}
\caption{\label{fig:trends_main}Increasing trends in the results reported in Table~\ref{tab:main_results}. Logarithms are base-10.}
\end{figure}

Figure~\ref{fig:trends_main} shows the main trends in the results. The optimality gap and running time both increase with the number of buses in the instances. Within each instance, there is also an upward trend as the number of periods increases from 4 to 24. Regarding the number of iterations, there also seems to be an increasing trend as a function of the number of periods for Case-6 instances, but not at all for the Case-24 ones. This is just because many of the associated executions reach the time limit and, thus, since the larger instances have more difficult SDP relaxations, the number of iterations that can be performed within the time limit goes down.

\subsection{Enhancing the analysis for the Case-118 instances}\label{sec:case118finer}

The goal of this section is to get a better understanding of the potential of the approach we have developed in this paper to tackle moderate-size instances such as Case-118. To this end, we concentrate on the most difficult variant, Case-118b, and define a series of variations of it as follows. We take as baseline the $P^{\text{min}}_i$ and $P^{\text{max}}_i$ values from \shortciteA{castillo2016}. Then, for each $\gamma \in \{ 1, 1.2, 1.4, 1.6, 1.8, 2, 3\}$, we replace these minimum and maximum values with $P^{\text{min}}_i/\gamma$ and $\gamma\cdot P^{\text{max}}_i$, respectively. Value $\gamma=1$ then leads to the values in \shortciteA{castillo2016} and, as $\gamma$ increases, we get closer to the MATPOWER \shortcite{zimmerman2010matpower} values, with which the Case-118 instances were solved at the root node. 

\begin{table}[!htbp]
\centering
{\small
\begin{tabular}[t]{lcccrcrc}
\toprule
\multicolumn{1}{c}{Periods} & \multicolumn{1}{c}{$\gamma$} & \multicolumn{1}{c}{MISDP-Gap} & \multicolumn{1}{c}{MIQCQP-Gap} & \multicolumn{1}{c}{Time} & \multicolumn{1}{c}{LS\,\%} & \multicolumn{1}{c}{Iter} & \multicolumn{1}{c}{UB jumps}\\
\midrule
\cellcolor{gray!25}{4} & \cellcolor{gray!25}{3.0} & \cellcolor{gray!25}{0.34\,\%} & \cellcolor{gray!25}{0.34\,\%} & \cellcolor{gray!25}{14.91} & \cellcolor{gray!25}{82.41\,\%} & \cellcolor{gray!25}{1.0} & \cellcolor{gray!25}{1.0}\\
4 & 2.0 & 0.84\,\% & 0.84\,\% & 24.65 & 79.03\,\% & 3.1 & 1.4\\
\cellcolor{gray!25}{4} & \cellcolor{gray!25}{1.8} & \cellcolor{gray!25}{1.68\,\%} & \cellcolor{gray!25}{1.68\,\%} & \cellcolor{gray!25}{19.60} & \cellcolor{gray!25}{83.54\,\%} & \cellcolor{gray!25}{1.3} & \cellcolor{gray!25}{1.1}\\
4 & 1.6 & 1.96\,\% & 1.96\,\% & 77.30 & 56.02\,\% & 14.7 & 1.4\\
\cellcolor{gray!25}{4} & \cellcolor{gray!25}{1.4} & \cellcolor{gray!25}{2.00\,\%} & \cellcolor{gray!25}{2.00\,\%} & \cellcolor{gray!25}{1234.37} & \cellcolor{gray!25}{19.85\,\%} & \cellcolor{gray!25}{785.4} & \cellcolor{gray!25}{2.7}\\
4 & 1.2 & 2.77\,\% & 2.77\,\% & \textbf{---} & 18.71\,\% & 2308.7 & 4.9\\
\cellcolor{gray!25}{4} & \cellcolor{gray!25}{1.0} & \cellcolor{gray!25}{6.15\,\%} & \cellcolor{gray!25}{6.15\,\%} & \cellcolor{gray!25}{\textbf{---}} & \cellcolor{gray!25}{29.89\,\%} & \cellcolor{gray!25}{2022.7} & \cellcolor{gray!25}{6.7}\\
12 & 3.0 & 0.39\,\% & 0.39\,\% & 227.52 & 64.15\,\% & 21.1 & 2.1\\
\cellcolor{gray!25}{12} & \cellcolor{gray!25}{2.0} & \cellcolor{gray!25}{1.07\,\%} & \cellcolor{gray!25}{1.07\,\%} & \cellcolor{gray!25}{128.38} & \cellcolor{gray!25}{84.04\,\%} & \cellcolor{gray!25}{3.1} & \cellcolor{gray!25}{1.7}\\
12 & 1.8 & 1.22\,\% & 1.22\,\% & 301.03 & 81.45\,\% & 10.0 & 2.7\\
\cellcolor{gray!25}{12} & \cellcolor{gray!25}{1.6} & \cellcolor{gray!25}{1.91\,\%} & \cellcolor{gray!25}{1.91\,\%} & \cellcolor{gray!25}{1957.60} & \cellcolor{gray!25}{78.11\,\%} & \cellcolor{gray!25}{77.9} & \cellcolor{gray!25}{3.3}\\
12 & 1.4 & 2.85\,\% & 2.85\,\% & \textbf{---} & 47.84\,\% & 380.6 & 3.1\\
\cellcolor{gray!25}{12} & \cellcolor{gray!25}{1.2} & \cellcolor{gray!25}{15.38\,\%} & \cellcolor{gray!25}{15.42\,\%} & \cellcolor{gray!25}{\textbf{---}} & \cellcolor{gray!25}{67.34\,\%} & \cellcolor{gray!25}{252.1} & \cellcolor{gray!25}{4.0}\\
12 & 1.0 & 23.28\,\% & 23.28\,\% & \textbf{---} & 79.65\,\% & 164.9 & 2.6\\
\cellcolor{gray!25}{24} & \cellcolor{gray!25}{3.0} & \cellcolor{gray!25}{1.13\,\%} & \cellcolor{gray!25}{1.13\,\%} & \cellcolor{gray!25}{228.86} & \cellcolor{gray!25}{88.20\,\%} & \cellcolor{gray!25}{1.7} & \cellcolor{gray!25}{1.3}\\
24 & 2.0 & 1.41\,\% & 1.41\,\% & 218.95 & 89.74\,\% & 1.3 & 1.1\\
\cellcolor{gray!25}{24} & \cellcolor{gray!25}{1.8} & \cellcolor{gray!25}{1.19\,\%} & \cellcolor{gray!25}{1.19\,\%} & \cellcolor{gray!25}{627.37} & \cellcolor{gray!25}{84.05\,\%} & \cellcolor{gray!25}{6.4} & \cellcolor{gray!25}{2.1}\\
24 & 1.6 & 1.81\,\% & 1.81\,\% & 1936.17 & 80.83\,\% & 25.0 & 3.3\\
\cellcolor{gray!25}{24} & \cellcolor{gray!25}{1.4} & \cellcolor{gray!25}{2.70\,\%} & \cellcolor{gray!25}{2.70\,\%} & \cellcolor{gray!25}{3321.95} & \cellcolor{gray!25}{68.64\,\%} & \cellcolor{gray!25}{76.3} & \cellcolor{gray!25}{2.7}\\
24 & 1.2 & 41.06\,\% & 41.06\,\% & \textbf{---} & 76.87\,\% & 69.0 & 1.1\\
\cellcolor{gray!25}{24} & \cellcolor{gray!25}{1.0} & \cellcolor{gray!25}{$\infty$} & \cellcolor{gray!25}{$\infty$} & \cellcolor{gray!25}{\textbf{---}} & \cellcolor{gray!25}{86.65\,\%} & \cellcolor{gray!25}{43.1} & \cellcolor{gray!25}{0.9}\\
\bottomrule
\end{tabular}
}
\caption{\label{tab:118_results} Average results by instance and number of periods in the Case-118b instances.}
\end{table}

Table~\ref{tab:118_results} reports the results of the computational analysis for these variants of the Case-118b instances. Again, for each number of periods and each value of $\gamma$, we have 7 variations of each instance run with 4 values of \texttt{run\_local}, and the best-performing execution among these 4 is taken to compute the averages shown in Table~\ref{tab:118_results}.

We can see that, as expected, the instances become more difficult as $\gamma$ goes down from~3 to~1, with increasing trends in all the performance metrics: optimality gaps, running times, iterations, and also UB jumps. Indeed, we can see that for values $1$, $1.2$, and $1.4$ no instance was solved, for all of its variations, with the specified tolerance within the time limit. Further, the $\infty$ values in the last row of Table~\ref{tab:118_results} are because, for 24 periods and $\gamma=1.0$, for some of the generated instances our algorithm could not find valid both lower and upper bounds (for any local solver configuration). 

\begin{table}[!htbp]
\centering
{\small
\begin{tabular}[t]{lccccrcrc}
\toprule
\multicolumn{1}{c}{Periods} & \multicolumn{1}{c}{Noise \%} & \multicolumn{1}{c}{Gen.} & \multicolumn{1}{c}{MISDP-Gap} & \multicolumn{1}{c}{MIQCQP-Gap} & \multicolumn{1}{c}{Time} & \multicolumn{1}{c}{LS\,\%} & \multicolumn{1}{c}{Iter} & \multicolumn{1}{c}{UB jumps}\\
\midrule
\cellcolor{gray!25}{24} & \cellcolor{gray!25}{0} & \cellcolor{gray!25}{\textbf{---}} & \cellcolor{gray!25}{39.98\,\%} & \cellcolor{gray!25}{39.98\,\%} & \cellcolor{gray!25}{\textbf{---}} & \cellcolor{gray!25}{96.82\,\%} & \cellcolor{gray!25}{9} & \cellcolor{gray!25}{2}\\
24 & 2 & 1 & 32.24\,\% & 32.24\,\% & \textbf{---} & 91.66\,\% & 25 & 2\\
\cellcolor{gray!25}{24} & \cellcolor{gray!25}{2} & \cellcolor{gray!25}{2} & \cellcolor{gray!25}{72.75\,\%} & \cellcolor{gray!25}{72.75\,\%} & \cellcolor{gray!25}{\textbf{---}} & \cellcolor{gray!25}{92.10\,\%} & \cellcolor{gray!25}{25} & \cellcolor{gray!25}{1}\\
24 & 4 & 1 & $\infty$ & $\infty$ & \textbf{---} & 75.93\,\% & 82 & 0\\
\cellcolor{gray!25}{24} & \cellcolor{gray!25}{4} & \cellcolor{gray!25}{2} & \cellcolor{gray!25}{93.68\,\%} & \cellcolor{gray!25}{93.68\,\%} & \cellcolor{gray!25}{\textbf{---}} & \cellcolor{gray!25}{69.01\,\%} & \cellcolor{gray!25}{101} & \cellcolor{gray!25}{0}\\
24 & 6 & 1 & 93.69\,\% & 93.69\,\% & \textbf{---} & 89.16\,\% & 35 & 1\\
\cellcolor{gray!25}{24} & \cellcolor{gray!25}{6} & \cellcolor{gray!25}{2} & \cellcolor{gray!25}{93.66\,\%} & \cellcolor{gray!25}{93.66\,\%} & \cellcolor{gray!25}{\textbf{---}} & \cellcolor{gray!25}{91.88\,\%} & \cellcolor{gray!25}{25} & \cellcolor{gray!25}{0}\\
\bottomrule
\end{tabular}
}
\caption{\label{tab:118_1.0_results}Disaggregate results for $\gamma=1$ and 24 periods.}
\end{table}

Table~\ref{tab:118_1.0_results} provides a disaggregated view of the best LS execution for each of the seven 24-period instances with $\gamma=1.0$. There was just one execution in which no gap was found, but the optimality gaps in the remaining cases are poor.\footnote{There are two executions where UB jumps is 0 and, yet, a finite gap is reported. This is because the feasible primal solution was not found during branch and bound, but in the initial call to the local solver at the root node.} Notably, in all instances a significant portion of time is spent in the local solver, confirming that the underlying MIQCQP problems remain challenging even if one gives up on certification of global optimality.

\begin{figure}[!htbp]
\centering
\includegraphics[width=\textwidth]{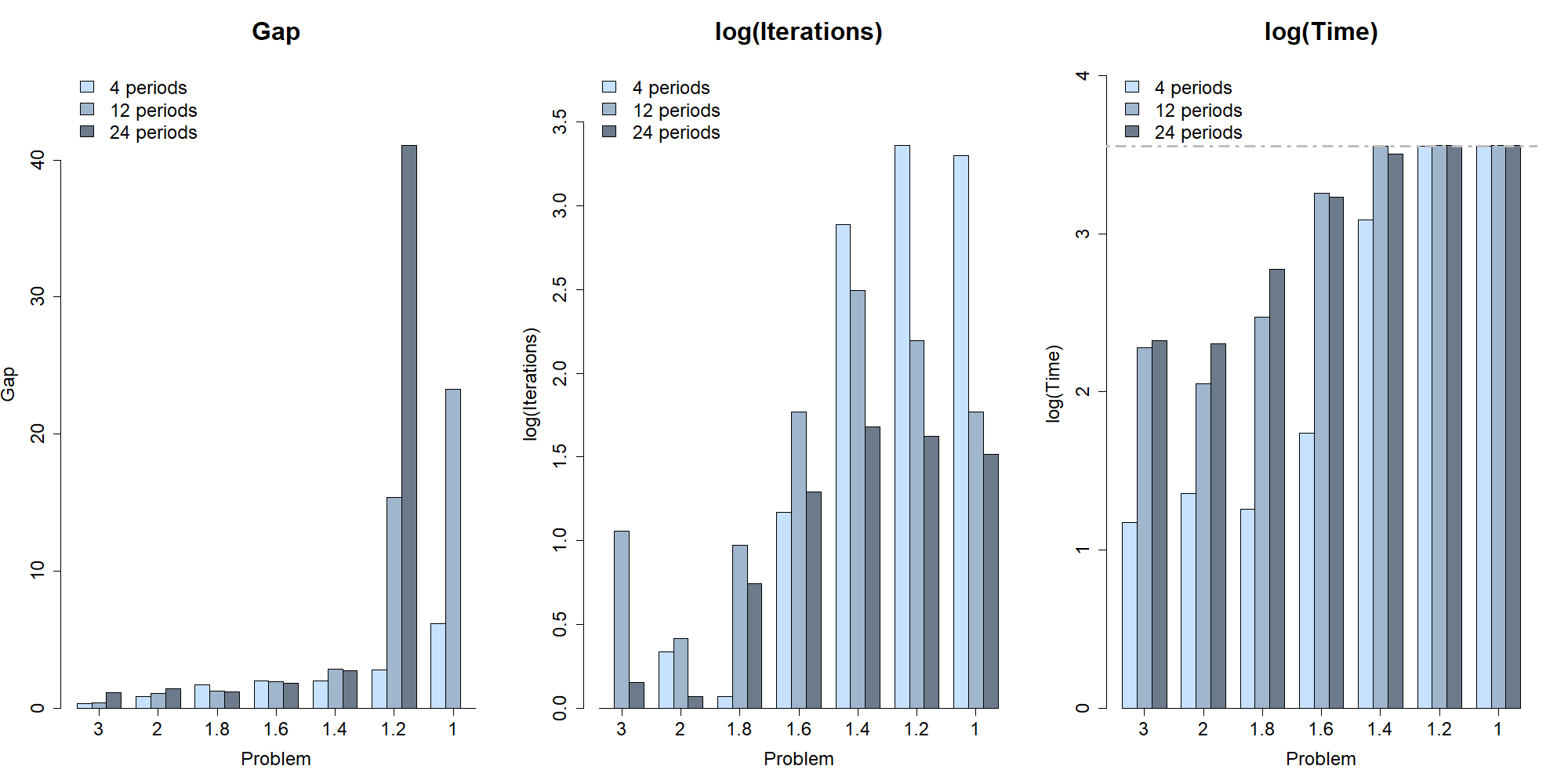}
\caption{\label{fig:trends_118}Increasing trends in the results reported in Table~\ref{tab:118_results}. Logarithms are base-10.}
\end{figure}

Figure~\ref{fig:trends_118} again shows the increasing trends from easier to harder instances. The effect of parameter $\gamma$ is evident in all subfigures, particularly in the optimality gaps.

Overall, the results in Table~\ref{tab:118_results} show that, despite the difficulties faced for the challenging instances resulting from low values of $\gamma$, our approach is quite effective at tackling Case-118 instances. As we observed in the previous tables, the quality of the [SDP-R] relaxations results in the gaps for problems [MISDP] and [MIQCQP] being close to each other.

\subsection{Local solver trade-offs}\label{sec:local}

In this section we present disaggregate results for the different values of \texttt{run\_local} parameter, to see how the frequency of the calls to the local solver impacts overall performance. The results show substantial variability, with no clear strategy to anticipate the best parameter choice for a given problem. This is why the results in the preceding sections already report the best configuration (ex post), as our main goal was to illustrate the potential of our approach.

\begin{table}[!htbp]
\centering
{\small
\begin{tabular}[t]{lcccrcrc}
\toprule
\multicolumn{1}{c}{Periods} & \multicolumn{1}{c}{\texttt{run\_local}} & \multicolumn{1}{c}{MISDP-Gap} & \multicolumn{1}{c}{MIQCQP-Gap} & \multicolumn{1}{c}{Time} & \multicolumn{1}{c}{LS\,\%} & \multicolumn{1}{c}{Iter} & \multicolumn{1}{c}{UB jumps}\\
\midrule
\cellcolor{gray!25}{4} & \cellcolor{gray!25}\text{best} & \cellcolor{gray!25}{1.38\,\%} & \cellcolor{gray!25}{2.14\,\%} & \cellcolor{gray!25}{9.44} & \cellcolor{gray!25}{61.52\,\%} & \cellcolor{gray!25}{42.75} & \cellcolor{gray!25}{2.15}\\
4 & 1.10 & 1.41\,\% & 2.17\,\% & 14.15 & 75.83\,\% & 42.30 & 2.50\\
4 & 1.25 & 1.40\,\% & 2.49\,\% & 8.11 & 65.71\,\% & 42.60 & 2.20\\
4 & 1.50 & 1.40\,\% & 3.11\,\% & 5.11 & 50.97\,\% & 42.90 & 1.70\\
4 & 2.00 & 1.37\,\% & $\infty$ & 3.90 & 38.85\,\% & 43.10 & 0.80\\
\cellcolor{gray!25}{12} & \cellcolor{gray!25}\text{best} & \cellcolor{gray!25}{1.37\,\%} & \cellcolor{gray!25}{2.03\,\%} & \cellcolor{gray!25}{97.67} & \cellcolor{gray!25}{22.16\,\%} & \cellcolor{gray!25}{481.15} & \cellcolor{gray!25}{2.65}\\
12 & 1.10 & 1.73\,\% & 2.04\,\% & 144.86 & 50.90\,\% & 454.60 & 3.60\\
12 & 1.25 & 1.94\,\% & 2.01\,\% & 105.76 & 32.35\,\% & 465.60 & 3.60\\
12 & 1.50 & 1.49\,\% & 2.06\,\% & 90.05 & 19.87\,\% & 471.70 & 2.70\\
12 & 2.00 & 1.39\,\% & $\infty$ & 84.91 & 9.89\,\% & 492.90 & 1.40\\
\cellcolor{gray!25}{24} & \cellcolor{gray!25}\text{best} & \cellcolor{gray!25}{1.94\,\%} & \cellcolor{gray!25}{1.94\,\%} & \cellcolor{gray!25}{311.40} & \cellcolor{gray!25}{23.42\,\%} & \cellcolor{gray!25}{898.85} & \cellcolor{gray!25}{3.00}\\
24 & 1.10 & 1.96\,\% & 1.96\,\% & 433.84 & 47.43\,\% & 884.70 & 4.00\\
24 & 1.25 & 1.95\,\% & 1.95\,\% & 322.60 & 27.72\,\% & 890.10 & 3.10\\
24 & 1.50 & 1.95\,\% & 1.95\,\% & 282.44 & 15.82\,\% & 893.70 & 2.90\\
24 & 2.00 & 1.98\,\% & 1.98\,\% & 265.01 & 6.21\,\% & 906.60 & 1.10\\
\bottomrule
\end{tabular}
}
\caption{\label{tab:LS6_results}Average results for Case-6a and Case-6b (averaging also over 6a and 6b).}
\end{table}

In Table~\ref{tab:LS6_results} we present the results for Case-6 instances. Recall that larger values of \texttt{run\_local} correspond to less frequent calls to the local solver and that the best LS configuration corresponds to the choice of \texttt{run\_local} with a smaller average between the gap for [MISDP] and the gap for [MIQCQP]. This implies that, by construction, the row of the best LS configuration has to be better than any other row at least in one of the two columns, but not necessarily in both of them.

As expected, the share of time spent by the local solver decreases as \texttt{run\_local} increases, while the number of UB jumps also declines. In terms of optimality gaps, performance is similar for all configurations, with lower values of \texttt{run\_local} producing slightly better gaps in general. When looking at running times, larger values of \texttt{run\_local} seem to be faster but, at the same time, there are some instances where they do not get feasible solutions for [MIQCQP], so the local solver calls from the [MISDP] solutions were not successful in these instances.

Note that, for the instances with 4 periods, the best LS Configuration delivers an average gap of 1.38, whereas the configuration with value 2 has an slightly better value, 1.37. As explained above, this arises because the best LS configuration is selected based on the average of MISDP-Gap and MIQCQP-Gap. Configuration 2 has been quite effective at finding good gaps for [MISDP] on some instances (more iterations due to fewer local solver calls), but yields relatively poor gaps for [MIQCQP] (fewer local solver calls).

\begin{table}[!htbp]
\centering
{\small
\begin{tabular}[t]{lcccrcrc}
\toprule
\multicolumn{1}{c}{Periods} & \multicolumn{1}{c}{\texttt{run\_local}} & \multicolumn{1}{c}{MISDP-Gap} & \multicolumn{1}{c}{MIQCQP-Gap} & \multicolumn{1}{c}{Time} & \multicolumn{1}{c}{LS\,\%} & \multicolumn{1}{c}{Iter} & \multicolumn{1}{c}{UB jumps}\\
\midrule
\cellcolor{gray!25}{4} & \cellcolor{gray!25}\text{best} & \cellcolor{gray!25}{4.62\,\%} & \cellcolor{gray!25}{4.63\,\%} & \cellcolor{gray!25}{2331.98} & \cellcolor{gray!25}{17.78\,\%} & \cellcolor{gray!25}{4236.7} & \cellcolor{gray!25}{6.35}\\
4 & 1.10 & 6.40\,\% & 6.41\,\% & 2341.80 & 26.71\,\% & 3892.5 & 7.80\\
4 & 1.25 & 11.58\,\% & 11.59\,\% & 2334.75 & 15.77\,\% & 4448.3 & 6.20\\
4 & 1.50 & 11.90\,\% & 11.90\,\% & 2544.21 & 5.96\,\% & 5137.6 & 5.40\\
4 & 2.00 & 13.04\,\% & 13.06\,\% & 2746.72 & 3.16\,\% & 5686.2 & 4.70\\
\cellcolor{gray!25}{12} & \cellcolor{gray!25}\text{best} & \cellcolor{gray!25}{2.90\,\%} & \cellcolor{gray!25}{2.90\,\%} & \cellcolor{gray!25}{2463.63} & \cellcolor{gray!25}{62.12\,\%} & \cellcolor{gray!25}{1190.3} & \cellcolor{gray!25}{5.00}\\
12 & 1.10 & 5.65\,\% & 5.65\,\% & 2624.12 & 95.28\,\% & 88.1 & 5.60\\
12 & 1.25 & 5.58\,\% & 5.58\,\% & 2494.62 & 70.91\,\% & 965.9 & 5.00\\
12 & 1.50 & 4.78\,\% & 4.78\,\% & 2410.97 & 59.44\,\% & 1283.4 & 4.60\\
12 & 2.00 & 4.18\,\% & 4.18\,\% & 2389.46 & 37.62\,\% & 1822.2 & 2.90\\
\cellcolor{gray!25}{24} & \cellcolor{gray!25}\text{best} & \cellcolor{gray!25}{6.39\,\%} & \cellcolor{gray!25}{6.39\,\%} & \cellcolor{gray!25}{3456.13} & \cellcolor{gray!25}{77.52\,\%} & \cellcolor{gray!25}{407.0} & \cellcolor{gray!25}{3.15}\\
24 & 1.10 & 9.39\,\% & 9.39\,\% & \textbf{---} & 98.11\,\% & 22.7 & 3.60\\
24 & 1.25 & 10.42\,\% & 10.42\,\% & 3577.77 & 92.07\,\% & 142.3 & 3.10\\
24 & 1.50 & 12.16\,\% & 12.16\,\% & \textbf{---} & 75.63\,\% & 460.2 & 2.90\\
24 & 2.00 & 13.95\,\% & 13.95\,\% & 3454.29 & 48.59\,\% & 955.5 & 2.40\\
\bottomrule
\end{tabular}
}
\caption{\label{tab:LS24_results}Average results for Case-24a and Case-24b (averaging also over 24a and 24b).}
\end{table}

Table~\ref{tab:LS24_results} shows the results for the 24 bus instances. Now, being able to choose the best-performing \texttt{run\_local} option for each of the 7 variants of each baseline instance significantly outperforms the performance obtained when sticking to the same value of \texttt{run\_local} for all 7 executions. The difference in performance is particularly noticeable in the optimality gaps. On the other hand, if we compared among them the different choices for \texttt{run\_local}, we observe that value $1.1$ seems to be the best choice for 4 and 24 periods, whereas $2$ is the best choice for 12 periods.

\begin{table}[!htbp]
\centering
{\small
\begin{tabular}[t]{lccccrcc}
\toprule
\multicolumn{1}{c}{Periods} & \multicolumn{1}{c}{\texttt{run\_local}} & \multicolumn{1}{c}{MISDP-Gap} & \multicolumn{1}{c}{MIQCQP-Gap} & \multicolumn{1}{c}{Time} & \multicolumn{1}{c}{LS\,\%} & \multicolumn{1}{c}{Iter} & \multicolumn{1}{c}{UB jumps}\\
\midrule
\cellcolor{gray!25}{4} & \cellcolor{gray!25}\text{best} & \cellcolor{gray!25}{1.36\,\%} & \cellcolor{gray!25}{1.36\,\%} & \cellcolor{gray!25}{274.17} & \cellcolor{gray!25}{64.17\,\%} & \cellcolor{gray!25}{161.10} & \cellcolor{gray!25}{1.52}\\
4 & 1.10 & 1.44\,\% & 1.44\,\% & 443.46 & 77.76\,\% & 135.70 & 1.60\\
4 & 1.25 & 1.37\,\% & 1.37\,\% & 298.28 & 73.07\,\% & 135.80 & 1.80\\
4 & 1.50 & 1.38\,\% & 1.38\,\% & 271.35 & 66.82\,\% & 165.90 & 1.40\\
4 & 2.00 & 1.39\,\% & 1.39\,\% & 317.55 & 57.59\,\% & 225.30 & 1.20\\
\cellcolor{gray!25}{12} & \cellcolor{gray!25}\text{best} & \cellcolor{gray!25}{1.49\,\%} & \cellcolor{gray!25}{1.49\,\%} & \cellcolor{gray!25}{1244.86} & \cellcolor{gray!25}{71.12\,\%} & \cellcolor{gray!25}{98.54} & \cellcolor{gray!25}{2.58}\\
12 & 1.10 & 1.81\,\% & 1.81\,\% & 1390.95 & 87.67\,\% & 23.40 & 2.70\\
12 & 1.25 & 1.62\,\% & 1.62\,\% & 1193.64 & 77.26\,\% & 92.10 & 2.60\\
12 & 1.50 & 1.73\,\% & 1.73\,\% & 1455.02 & 60.10\,\% & 212.70 & 2.70\\
12 & 2.00 & 1.89\,\% & 1.89\,\% & 1557.61 & 45.01\,\% & 258.70 & 2.20\\
\cellcolor{gray!25}{24} & \cellcolor{gray!25}\text{best} & \cellcolor{gray!25}{1.65\,\%} & \cellcolor{gray!25}{1.65\,\%} & \cellcolor{gray!25}{1266.66} & \cellcolor{gray!25}{82.29\,\%} & \cellcolor{gray!25}{22.14} & \cellcolor{gray!25}{2.10}\\
24 & 1.10 & 1.96\,\% & 1.91\,\% & 1691.22 & 91.09\,\% & 7.40 & 2.30\\
24 & 1.25 & 1.97\,\% & 1.97\,\% & 1488.11 & 86.49\,\% & 17.40 & 2.10\\
24 & 1.50 & 1.82\,\% & 1.87\,\% & 1304.25 & 78.57\,\% & 31.40 & 2.00\\
24 & 2.00 & 2.04\,\% & 2.04\,\% & 1579.71 & 58.03\,\% & 81.90 & 1.90\\
\bottomrule
\end{tabular}
}
\caption{\label{tab:LS118_results}Average results for Case-118b (removing values 1 and 1.2 for $\gamma$).}
\end{table}

Finally, Table~\ref{tab:LS118_results} shows the results for Case-118b instances, averaging over the different values of $\gamma$. Yet, in order to get comparable results for the different numbers of periods and LS configurations, we remove the instances with $\gamma=1$ since, as we have seen, in some instances with 24 periods they do not get a finite gap within the time limit. Also, we remove instances with $\gamma=1.2$, since no LS configuration was able to solve all of them either (each instance with $\gamma=1.2$ was solved by at least one LS configuration, but no configuration solved all of them). By removing the most difficult instances we can see more clearly the behavior with respect to the LS configuration.

Again, the best configuration outperforms the fixed-\texttt{run\_local} executions in terms of gap, although not as clearly as for Case-24. Also, in terms of running time, the best configuration does not beat all fixed-\texttt{run\_local} configurations. Regarding these configurations, $1.25$ and $1.5$ seem to be the best choices for 4~periods, $1.25$ for 12~periods and $1.5$ for 24~periods. Interestingly, this situation is quite opposite from what we had for Case-6 and Case-24, where the best fixed-\texttt{run\_local} configurations were precisely $1.1$ and $2$.

Summing up, we have seen that the choice of the frequency of calls to the local solver can have a significant impact on performance, but that fine-tuning this parameter may not be an easy task.

\subsection{Comparison with state-of-the-art solvers}

Given the encouraging performance observed thus far, we now assess the extent to which the combination of the different methodologies behind Algorithm~\ref{algorithm:full} is competitive with a direct solution with leading commercial solvers. To this end, in this section we revisit the average results presented in Table~\ref{tab:main_results} and compare them with those obtained when tackling the same instances with Gurobi~13.0.0 \shortcite{gurobi2025} and BARON~25.12.10 \shortcite{baron2025}.

\begin{table}[!htbp]
\centering
{\small
\begin{tabular}[t]{lccccrrr}
\toprule
\multicolumn{1}{c}{Instance} & \multicolumn{1}{c}{Periods} & \multicolumn{1}{c}{MIQCQP-Gap} & \multicolumn{1}{c}{Gap-GUR} & \multicolumn{1}{c}{Gap-BAR} & \multicolumn{1}{c}{Time} & \multicolumn{1}{c}{Time-GUR} & \multicolumn{1}{c}{Time-BAR}\\
\midrule
\cellcolor{gray!25}{Case-6a} & \cellcolor{gray!25}{4} & \cellcolor{gray!25}{1.17\,\%} & \cellcolor{gray!25}{25.05\,\%} & \cellcolor{gray!25}{0.81\,\%} & \cellcolor{gray!25}{3.60} & \cellcolor{gray!25}{2062.46} & \cellcolor{gray!25}{1038.38}\\
Case-6a & 12 & 2.07\,\% & 2.54\,\% & 4.20\,\% & 43.30 & 2127.67 & \textbf{---}\\
\cellcolor{gray!25}{Case-6a} & \cellcolor{gray!25}{24} & \cellcolor{gray!25}{1.90\,\%} & \cellcolor{gray!25}{$\infty$} & \cellcolor{gray!25}{$\infty$} & \cellcolor{gray!25}{173.74} & \cellcolor{gray!25}{852.82} & \cellcolor{gray!25}{\textbf{---}}\\
Case-6b & 4 & 3.11\,\% & 1.89\,\% & 6.76\,\% & 15.27 & 794.69 & \textbf{---}\\
\cellcolor{gray!25}{Case-6b} & \cellcolor{gray!25}{12} & \cellcolor{gray!25}{1.99\,\%} & \cellcolor{gray!25}{1.78\,\%} & \cellcolor{gray!25}{7.56\,\%} & \cellcolor{gray!25}{152.04} & \cellcolor{gray!25}{285.85} & \cellcolor{gray!25}{\textbf{---}}\\
Case-6b & 24 & 1.97\,\% & 1.73\,\% &$ \infty$ & 449.06 & 489.17 & \textbf{---}\\
\cellcolor{gray!25}{Case-24a} & \cellcolor{gray!25}{4} & \cellcolor{gray!25}{1.84\,\%} & \cellcolor{gray!25}{6.97\,\%} & \cellcolor{gray!25}{$\infty$} & \cellcolor{gray!25}{1063.67} & \cellcolor{gray!25}{\textbf{---}} & \cellcolor{gray!25}{\textbf{---}}\\
Case-24a & 12 & 1.88\,\% & $\infty$ & $\infty$ & 1326.12 & \textbf{---} & \textbf{---}\\
\cellcolor{gray!25}{Case-24a} & \cellcolor{gray!25}{24} & \cellcolor{gray!25}{3.33\,\%} & \cellcolor{gray!25}{$\infty$} & \cellcolor{gray!25}{$\infty$} & \cellcolor{gray!25}{3309.43} & \cellcolor{gray!25}{\textbf{---}} & \cellcolor{gray!25}{\textbf{---}}\\
Case-24b & 4 & 7.41\,\% & 18.40\,\% & $\infty$ & \textbf{---} & \textbf{---} & \textbf{---}\\
\cellcolor{gray!25}{Case-24b} & \cellcolor{gray!25}{12} & \cellcolor{gray!25}{3.91\,\%} & \cellcolor{gray!25}{$\infty$} & \cellcolor{gray!25}{$\infty$} & \cellcolor{gray!25}{\textbf{---}} & \cellcolor{gray!25}{\textbf{---}} & \cellcolor{gray!25}{\textbf{---}}\\
Case-24b & 24 & 9.45\,\% & $\infty$ & $\infty$ & \textbf{---} & \textbf{---} & \textbf{---}\\
\cellcolor{gray!25}{Case-118a} & \cellcolor{gray!25}{4} & \cellcolor{gray!25}{0.01\,\%} & \cellcolor{gray!25}{92.50\,\%} & \cellcolor{gray!25}{$\infty$} & \cellcolor{gray!25}{17.40} & \cellcolor{gray!25}{\textbf{---}} & \cellcolor{gray!25}{\textbf{---}}\\
Case-118a & 12 & 0.06\,\% & 98.06\,\% & $\infty$ & 56.84 & \textbf{---} & \textbf{---}\\
\cellcolor{gray!25}{Case-118a} & \cellcolor{gray!25}{24} & \cellcolor{gray!25}{0.04\,\%} & \cellcolor{gray!25}{99.19\,\%} & \cellcolor{gray!25}{$\infty$} & \cellcolor{gray!25}{137.84} & \cellcolor{gray!25}{\textbf{---}} & \cellcolor{gray!25}{\textbf{---}}\\
Case-118b & 4 & 0.00\,\% & 100.00\,\% & $\infty$ & 13.32 & \textbf{---} & \textbf{---}\\
\cellcolor{gray!25}{Case-118b} & \cellcolor{gray!25}{12} & \cellcolor{gray!25}{0.02\,\%} & \cellcolor{gray!25}{100.00\,\%} & \cellcolor{gray!25}{$\infty$} & \cellcolor{gray!25}{45.58} & \cellcolor{gray!25}{\textbf{---}} & \cellcolor{gray!25}{\textbf{---}}\\
Case-118b & 24 & 0.09\,\% & 100.00\,\% & $\infty$ & 126.68 & \textbf{---} & \textbf{---}\\
\bottomrule
\end{tabular}
}
\caption{\label{tab:solvers}Comparison of the average results obtained with Algorithm~1 (MIQCQP-Gap and Time) with those obtained with state-of-the-art global optimization solvers Gurobi~13.0.0 and BARON~25.12.10.}
\end{table}

The joint results are presented in Table~\ref{tab:solvers}. They show that Algorithm~\ref{algorithm:full} outperforms both Gurobi and BARON for most of the instances considered. The only exception is Case-6b, where Gurobi achieves the best relative gaps, although at the expense of longer running times. In all other instances, Algorithm~\ref{algorithm:full} clearly performs better, both in terms of optimality gaps and computational times.

A quick inspection of the running times already indicates that one of the main goals of this paper has been achieved: designing an algorithm that scales well as the size of the instances increases. Indeed, given the poor results obtained by both Gurobi and BARON on Case-118a and Case-118b instances, we did not test them on the more challenging variants introduced in Section~\ref{sec:case118finer}.

\section{Conclusions}\label{sec:conclusion}

Our primary focus in this work has been finding optimal solutions to MIQCQPs, under the general assumption that the quadratic objective and constraints are nonconvex, \emph{i.e.}, not all quadratic forms involved have positive semidefinite Hessian. This class of problems combines two types of nonconvexity: discreteness of a subset of variables and nonconvex nonlinear constraints and/or objective function. While branch-and-bound is the algorithm of choice to find optimal solutions of nonconvex optimization problems, plain branch-and-bound implementations that use simple McCormick relaxations of the quadratic constraints show poor performance due to the loose relaxation that only considers single bilinear terms. The SDP relaxation of a nonconvex QCQP provides much tighter bounds, but suffers from scalability problems and the lack of warm-starting mechanisms, the latter preventing its effective use in a branch-and-bound framework.

We propose a solution algorithm that exploits a key fact about several practical applications: their MIQCQP formulation has a sparse quadratic structure. Specifically, we use term and correlative sparsity in these quadratic forms to reformulate a large-scale SDP relaxation such that it can be solved more efficiently by enforcing positive semidefiniteness of much smaller submatrices via chordal extensions. Our solution approach embeds such sparse SDP solver, which uses the chordal extension tool CS-TSSOS, within a branch-and-bound algorithm, thus obtaining tight bounds and faster node solves. We add calls to an external local solver that helps find feasible solutions.

We test our approach on one of the toughest classes of MIQCQP, namely the Unit Commitment problem under ACOPF constraints. While Unit Commitment introduces binary variables and logical constraints for the operation of a power grid across multiple time periods, ACOPF constraints are quadratic, nonconvex, and sparse due to the sparse topology of the power grid. This is perhaps a quintessential hard MIQCQP in that it combines two problems that are NP-hard in their own right: optimization under ACOPF constraints and Unit Commitment. Our approach can handle problems with up to 118 nodes, 186 edges, and 24 time periods while maintaining a small optimality gap. Importantly, the proposed algorithm outperforms state-of-the-art global optimization solvers, with the performance advantage becoming more pronounced as problem sizes increase, which highlights the superior scaling behavior of our approach.

These results attest to the effectiveness of our proposed algorithm. Due to the complexity of some of the different methodologies integrated in the proposed approach, in our discussion of the computational results we provide a detailed account of the influence of its components: tightness of the SDP relaxation, performance when using correlative/term sparsity, frequency of the calls to a local solver, as well as the sensitivity of our solver to variants of the original instances.

The results we have obtained suggest possible research directions for this method. First, our approach is based on a simple branch-and-bound implementation that does not incorporate all routines that are critical to its effectiveness, such as presolving and cutting. Moreover, it uses an external SDP solver (although very efficient) and a separate tool to handle term and correlative sparsity. It would be interesting to know whether a more tightly integrated solver, based on a state-of-the-art branch-and-cut solver, can further improve the performance on nonconvex MIQCQP.

Another possible research direction consists of replacing the full resolution of the reduced-size SDP relaxation by adding SDP cuts  \shortcite{sherali2002enhancing} based on the same decomposition. This would take advantage of warm-starting capabilities in solving SDP  relaxations that are crucial in branch-and-cut algorithms. Finally, while our approach restricts branching rules to integer variables, it would be interesting to extend branching to continuous variables that appear in bilinear terms, similar to {\em spatial} branching mechanisms that are typical of nonconvex MINLO solvers.

\section*{Acknowledgments}

\footnotesize{{This work is part of the R\&D project PID2021-124030NB-C32 granted by MICIU/AEI/10.13039/501100011033. This research was also funded by Grupos de Referencia Competitiva ED431C-2021/24 from the Consellería de Cultura, Educación e Universidades, Xunta de Galicia. Ignacio Gómez-Casares acknowledges the support from the Spanish Ministry of Education through FPU grant 20/01555. Bissan Ghaddar's research is supported by the Natural Sciences and Engineering Research Council of Canada Discovery RGPIN-2025-04585 and by the John Thompson Chair Fellowship.
Pietro Belotti's work is supported by the HEXAGON project (CUP F53D23010010001), part of the National Recovery and Resilience Plan (NRRP), Mission 4, Component 2, Investment 1.1, Call for tender No.\mbox{} 1409, funded by the European Union –- NextGenerationEU.}
\bibliographystyle{apacite} 
\bibliography{references}

@article{faria2011,
  title={Novel bound contraction procedure for global optimization of bilinear {MINLP} problems with applications to water management problems},
  author={Faria, D{\'e}bora C and Bagajewicz, Miguel J},
  journal={Computers \& Chemical Engineering},
  volume={35},
  number={3},
  pages={446--455},
  year={2011},
  publisher={Elsevier}
}

@article{misener2009,
  title={Advances for the pooling problem: Modeling, global optimization, and computational studies},
  author={Misener, Ruth and Floudas, Christodoulos A},
  journal={Applied and Computational Mathematics},
  volume={8},
  number={1},
  pages={3--22},
  year={2009}
}

@article{liers2021,
  title={Solving mixed-integer nonlinear optimization problems using simultaneous convexification: a case study for gas networks},
  author={Liers, Frauke and Martin, Alexander and Merkert, Maximilian and Mertens, Nick and Michaels, Dennis},
  journal={Journal of Global Optimization},
  volume={80},
  number={2},
  pages={307--340},
  year={2021},
  publisher={Springer}
}

@article{li2024,
  title={A reformulation-enumeration {MINLP} algorithm for gas network design},
  author={Li, Yijiang and Dey, Santanu S and Sahinidis, Nikolaos V},
  journal={Journal of Global Optimization},
  volume={90},
  number={4},
  pages={931--963},
  year={2024},
  publisher={Springer}
}

@article{bragalli2012,
  title={On the optimal design of water distribution networks: a practical {MINLP} approach},
  author={Bragalli, Cristiana and D’Ambrosio, Claudia and Lee, Jon and Lodi, Andrea and Toth, Paolo},
  journal={Optimization and Engineering},
  volume={13},
  pages={219--246},
  year={2012},
  publisher={Springer}
}

@article{kocuk2017,
  title={New formulation and strong {MISOCP} relaxations for {AC} optimal transmission switching problem},
  author={Kocuk, Burak and Dey, Santanu S and Sun, Xu Andy},
  journal={IEEE Transactions on Power Systems},
  volume={32},
  number={6},
  pages={4161--4170},
  year={2017},
  publisher={IEEE}
}

@article{kocuk2018,
  title={Matrix minor reformulation and {SOCP}-based spatial branch-and-cut method for the {AC} optimal power flow problem},
  author={Kocuk, Burak and Dey, Santanu S and Sun, Xu Andy},
  journal={Mathematical Programming Computation},
  volume={10},
  number={4},
  pages={557--596},
  year={2018},
  publisher={Springer}
}

@article{bingane2019,
  title={Tight-and-cheap conic relaxation for the optimal reactive power dispatch problem},
  author={Bingane, Christian and Anjos, Miguel F and Le Digabel, S{\'e}bastien},
  journal={IEEE Transactions on Power Systems},
  volume={34},
  number={6},
  pages={4684--4693},
  year={2019},
  publisher={IEEE}
}

@article{Lavaei:2012,
	author = {Lavaei, Javad and Low, Steven H.},
	journal = {IEEE Transactions on Power Systems},
	number = {1},
	pages = {92--107},
	title = {Zero Duality Gap in Optimal Power Flow Problem},
	volume = {27},
	year = {2012}}

@article{Carpentier:1962aa,
	author = {Carpentier, Jacques},
	journal = {Bulletin de la Société Francaise des électriciens},
	number = {3},
	pages = {431-447},
	title = {Contribution á l'étude du dispatching économique},
	volume = {8},
	year = {1962}}

@article{Lehmann:2016aa,
	author = {Lehmann, Karsten and Grastien, Alban and Van Hentenryck, Pascal},
	journal = {IEEE Transactions on Power Systems},
	number = {1},
	pages = {798--801},
	title = {{AC}-Feasibility on Tree Networks is {NP}-Hard},
	volume = {31},
	year = {2016}}

@article{Frank:2012aa,
	author = {Frank, Stephen and Steponavice, Ingrida and Rebennack, Steffen},
	journal = {Energy Systems},
	number = {3},
	pages = {221--258},
	title = {Optimal power flow: a bibliographic survey. {I}. {F}ormulations and deterministic methods},
	volume = {3},
	year = {2012}}

@article{Garver:1962,
	author = {Garver, Len L.},
	journal = {Transactions of the American Institute of Electrical Engineers. Part III: Power Apparatus and Systems},
	number = {3},
	pages = {730--734},
	title = {Power Generation Scheduling by Integer Programming-Development of Theory},
	volume = {81},
	year = {1962}}

@article{Montero:2022,
	author = {Montero, Luis and Bello, Antonio and Reneses, Javier},
	journal = {Energies},
	number = {4},
	pages = {1296},
	title = {A Review on the Unit Commitment Problem: Approaches, Techniques, and Resolution Methods},
	volume = {15},
	year = {2022}}

@article{Ostrowski:2012,
	author = {Ostrowski, James and Anjos, Miguel F. and Vannelli, Anthony},
	journal = {IEEE Transactions on Power Systems},
	number = {1},
	pages = {39--46},
	title = {Tight Mixed Integer Linear Programming Formulations for the Unit Commitment Problem},
	volume = {27},
	year = {2012}}

@article{Marcovecchio:2014,
	author = {Marcovecchio, Marian G. and Novais, Augusto Q. and Grossmann, Ignacio E.},
	journal = {Computers \& Chemical Engineering},
	pages = {53--68},
	title = {Deterministic optimization of the thermal Unit Commitment problem: A Branch and Cut search},
	volume = {67},
	year = {2014}}

@article{Damc-Kurt:2016,
	author = {Damcı-Kurt, Pelin and Küçükyavuz, Simge and Rajan, Deepak and Atamtürk, Alper},
	journal = {Mathematical Programming},
	number = {1-2},
	pages = {175--205},
	title = {A polyhedral study of production ramping},
	volume = {158},
	year = {2016}}

@article{Knueven:2018,
	author = {Knueven, Ben and Ostrowski, Jim and Wang, Jianhui},
	journal = {INFORMS Journal on Computing},
	number = {4},
	pages = {739--749},
	title = {The Ramping Polytope and Cut Generation for the Unit Commitment Problem},
	volume = {30},
	year = {2018}}

@article{Ostrowski:2015,
	author = {Ostrowski, James and Anjos, Miguel F. and Vannelli, Anthony},
	journal = {Mathematical Programming},
	number = {1},
	pages = {99--129},
	title = {Modified orbital branching for structured symmetry with an application to unit commitment},
	volume = {150},
	year = {2015}}

@article{Wang:2022,
	journal = {ACM Transactions on Mathematical Software},
	number = {4},
	pages = {1--26},
	title = {{CS-TSSOS}: Correlative and Term Sparsity for Large-Scale Polynomial Optimization},
    author={Wang, Jie and Magron, Victor and Lasserre, Jean-Bernard and Mai, Ngoc Hoang Anh},
	volume = {48},
	year = {2022}}

@article{chowdhury2025,
  title={Optimal Power Flow ({OPF}) Analysis for {AC}--{DC} Active Distribution Networks Utilizing Second-Order Cone Programming ({SOCP}) Approach},
  author={Chowdhury, Md Mahmud-Ul-Tarik and Murari, Krishna and Hasan, Md Shamim and Kamalasadan, Sukumar},
  journal={IEEE Transactions on Industrial Informatics},
  year={2025},
  volume={21},
  number={5},
  pages={3555--3564},
  publisher={IEEE}
}

@article{wang2020,
  title={Chordal-{TSSOS}: A {M}oment-{SOS} hierarchy that exploits term sparsity with chordal extension},
  author={Wang, Jie and Magron, Victor and Lasserre, Jean-Bernard},
  journal={SIAM Journal on Optimization},
  volume={31},
  number={1},
  pages={114--141},
  year={2021},
  publisher={SIAM}
}

@ARTICLE{Ghaddar2016, 
author={Bissan Ghaddar and Jakub Marecek and Martin Mevissen}, 
journal={{IEEE} Transactions on Power Systems}, 
title={Optimal Power Flow as a Polynomial Optimization Problem}, 
year={2016}, 
volume={31}, 
number={1}, 
pages={539-546}
}

@article{Ghaddar2017,
author = "Bissan Ghaddar and Mathieu Claeys and Martin Mevissen and Bradley J. Eck",
title = "Polynomial optimization for water networks: Global solutions for the valve setting problem",
journal = "European Journal of Operational Research",
volume = "261",
number = "2",
pages = "450 -- 459",
year = "2017"
}

@article{castillo2016,
  title={The unit commitment problem with {AC} optimal power flow constraints},
  author={Castillo, Anya and Laird, Carl and Silva-Monroy, C{\'e}sar A and Watson, Jean-Paul and O’Neill, Richard P},
  journal={IEEE Transactions on Power Systems},
  volume={31},
  number={6},
  pages={4853--4866},
  year={2016},
  publisher={IEEE}
}

@article{zimmerman2010matpower,
  title={{MATPOWER}: Steady-state operations, planning, and analysis tools for power systems research and education},
  author={Zimmerman, Ray Daniel and Murillo-S{\'a}nchez, Carlos Edmundo and Thomas, Robert John},
  journal={IEEE Transactions on Power Systems},
  volume={26},
  number={1},
  pages={12--19},
  year={2010},
  publisher={IEEE}
}

@book{anjos2011handbook,
  title={Handbook on semidefinite, conic and polynomial optimization},
  author={Anjos, Miguel F and Lasserre, Jean B},
  volume={166},
  year={2011},
  publisher={Springer Science \& Business Media}
}

@article{kim2003exact,
  title={Exact solutions of some nonconvex quadratic optimization problems via {SDP} and {SOCP} relaxations},
  author={Kim, Sunyoung and Kojima, Masakazu},
  journal={Computational Optimization and Applications},
  volume={26},
  pages={143--154},
  year={2003},
  publisher={Springer}
}

@article{waki2006sums,
  title={Sums of squares and semidefinite program relaxations for polynomial optimization problems with structured sparsity},
  author={Waki, Hayato and Kim, Sunyoung and Kojima, Masakazu and Muramatsu, Masakazu},
  journal={SIAM Journal on Optimization},
  volume={17},
  number={1},
  pages={218--242},
  year={2006},
  publisher={SIAM}
}

@article{bonami2018globally,
  title={Globally solving nonconvex quadratic programming problems with box constraints via integer programming methods},
  author={Bonami, Pierre and G{\"u}nl{\"u}k, Oktay and Linderoth, Jeff},
  journal={Mathematical Programming Computation},
  volume={10},
  number={3},
  pages={333--382},
  year={2018},
  publisher={Springer}
}

@article{saxena2010convex,
  title={Convex relaxations of non-convex mixed integer quadratically constrained programs: extended formulations},
  author={Saxena, Anureet and Bonami, Pierre and Lee, Jon},
  journal={Mathematical Programming},
  volume={124},
  number={1},
  pages={383--411},
  year={2010},
  publisher={Springer}
}

@article{saxena2011convex,
  title={Convex relaxations of non-convex mixed integer quadratically constrained programs: projected formulations},
  author={Saxena, Anureet and Bonami, Pierre and Lee, Jon},
  journal={Mathematical Programming},
  volume={130},
  number={2},
  pages={359--413},
  year={2011},
  publisher={Springer}
}

@article{anstreicher2009semidefinite,
  title={Semidefinite programming versus the reformulation-linearization technique for nonconvex quadratically constrained quadratic programming},
  author={Anstreicher, Kurt M},
  journal={Journal of Global Optimization},
  volume={43},
  number={2},
  pages={471--484},
  year={2009},
  publisher={Springer}
}

@article{adjiman1998global,
  title={A global optimization method, $\alpha${BB}, for general twice-differentiable constrained {NLP}s--{I}. {T}heoretical advances},
  author={Adjiman, Claire S and Dallwig, Stefan and Floudas, Christodoulos A and Neumaier, Arnold},
  journal={Computers \& Chemical Engineering},
  volume={22},
  number={9},
  pages={1137--1158},
  year={1998},
  publisher={Elsevier}
}

@article{bertsimas2021unified,
  title={A unified approach to mixed-integer optimization problems with logical constraints},
  author={Bertsimas, Dimitris and Cory-Wright, Ryan and Pauphilet, Jean},
  journal={SIAM Journal on Optimization},
  volume={31},
  number={3},
  pages={2340--2367},
  year={2021},
  publisher={SIAM}
}

@incollection{benson2013mixed,
  title={Mixed-integer second-order cone programming: A survey},
  author={Benson, Hande Y and Sa{\u{g}}lam, {\"U}mit},
  booktitle={Theory driven by influential applications},
  pages={13--36},
  year={2013},
  publisher={INFORMS}
}

@article{buchheim2013semidefinite,
  title={Semidefinite relaxations for non-convex quadratic mixed-integer programming},
  author={Buchheim, Christoph and Wiegele, Angelika},
  journal={Mathematical Programming},
  volume={141},
  number={1},
  pages={435--452},
  year={2013},
  publisher={Springer}
}

@article{lasserre2006convergent,
  title={Convergent {SDP}-relaxations in polynomial optimization with sparsity},
  author={Lasserre, Jean B},
  journal={SIAM Journal on Optimization},
  volume={17},
  number={3},
  pages={822--843},
  year={2006},
  publisher={SIAM}
}

@article{waki2008algorithm,
  title={Algorithm 883: Sparse{POP}---a sparse semidefinite programming relaxation of polynomial optimization problems},
  author={Waki, Hayato and Kim, Sunyoung and Kojima, Masakazu and Muramatsu, Masakazu and Sugimoto, Hiroshi},
  journal={ACM Transactions on Mathematical Software (TOMS)},
  volume={35},
  number={2},
  pages={1--13},
  year={2008},
  publisher={ACM New York, NY, USA}
}

@manual{mosek,
   author = "{MOSEK ApS}",
   year = 2019,
   note= {\url{https://www.mosek.com/}}
 }

@article{Carlson:2012aa,
	journal = {Interfaces},
	number = {1},
	pages = {58--73},
	title = {{MISO} Unlocks Billions in Savings Through the Application of Operations Research for Energy and Ancillary Services Markets},
	volume = {42},
    author = {Brian Carlson and  Yonghong Chen and Mingguo Hong and Roy Jones and Kevin Larson and Xingwang Ma and Peter Nieuwesteeg and Haili Song and Kimberly Sperry and Matthew Tackett and Doug Taylor and Jie Wan and Eugene Zak},
	year = {2012}
}

@article{Ma:1999aa,
	author = {Ma, Haili and Shahidehpour, S.M.},
	journal = {IEEE Transactions on Power Systems},
	number = {2},
	pages = {757--764},
	title = {Unit commitment with transmission security and voltage constraints},
	volume = {14},
	year = {1999}}

@inproceedings{Murillo-Sanchez:1999aa,
	author = {Murillo-Sanchez, C. and Thomas, R.J.},
	booktitle = {Thermal unit commitment with nonlinear power flow constraints},
	pages = {10},
	publisher = {IEEE},
	title = {Thermal unit commitment with nonlinear power flow constraints},
	year = {1999}}

@article{Nasri:2016aa,
	author = {Nasri, Amin and Kazempour, S. Jalal and Conejo, Antonio J. and Ghandhari, Mehrdad},
	journal = {IEEE Transactions on Power Systems},
	number = {1},
	pages = {412--422},
	title = {Network-Constrained {AC} Unit Commitment Under Uncertainty: A {B}enders' Decomposition Approach},
	volume = {31},
	year = {2016}}

@article{Sifuentes:2007aa,
	author = {Sifuentes, Wilfredo S. and Vargas, Alberto},
	journal = {IEEE Transactions on Power Systems},
	number = {3},
	pages = {1351--1359},
	title = {Hydrothermal Scheduling Using {B}enders Decomposition: Accelerating Techniques},
	volume = {22},
	year = {2007}}

@article{Fu:2006aa,
	author = {Fu, Yong and Shahidehpour, Mohammad and Li, Zuyi},
	journal = {IEEE Transactions on Power Systems},
	number = {2},
	pages = {897--908},
	title = {{AC} Contingency Dispatch Based on Security-Constrained Unit Commitment},
	volume = {21},
	year = {2006}}

@inproceedings{bhardwaj2012unit,
  title={Unit commitment in electrical power system -- a literature review},
  author={Bhardwaj, Amit and Kamboj, Vikram Kumar and Shukla, Vijay Kumar and Singh, Bhupinder and Khurana, Preeti},
  booktitle={2012 {IEEE} international power engineering and optimization conference -- {M}elaka, {M}alaysia},
  pages={275--280},
  year={2012},
  organization={IEEE}
}

@article{Padhy:2004aa,
	author = {Padhy, N.P.},
	journal = {IEEE Transactions on Power Systems},
	number = {2},
	pages = {1196--1205},
	title = {Unit Commitment---A Bibliographical Survey},
	volume = {19},
	year = {2004}}

@incollection{Byrd:2006,
	address = {Boston, MA},
	author = {Byrd, Richard H. and Nocedal, Jorge and Waltz, Richard A.},
	booktitle = {Large-Scale Nonlinear Optimization},
	doi = {10.1007/0-387-30065-1_4},
	editor = {Di Pillo, G. and Roma, M.},
	isbn = {978-0-387-30065-8},
	pages = {35--59},
	publisher = {Springer US},
	title = {Knitro: An Integrated Package for Nonlinear Optimization},
	year = {2006}}

@article{Bezanson:2017,
	author = {Bezanson, Jeff and Edelman, Alan and Karpinski, Stefan and Shah, Viral B},
	doi = {10.1137/141000671},
	journal = {SIAM Review},
	number = {1},
	pages = {65--98},
	publisher = {SIAM},
	title = {Julia: A fresh approach to numerical computing},
	volume = {59},
	year = {2017}}

@misc{python,
	author = {{Python Software Foundation}},
	note = {Available at: \url{http://www.python.org}},
	title = {Python Language Reference, version 3.11},
	year = {2024}}

@book{Bynum:2021aa,
	author = {Bynum, Michael L. and Hackebeil, Gabriel A. and Hart, William E. and Laird, Carl D. and Nicholson, Bethany L. and Siirola, John D. and Watson, Jean-Paul and Woodruff, David L.},
	edition = {Third},
	publisher = {Springer Science \& Business Media},
	title = {Pyomo--optimization modeling in {P}ython},
	volume = {67},
	year = {2021}}

@article{Hart:2011aa,
	author = {Hart, William E and Watson, Jean-Paul and Woodruff, David L},
	journal = {Mathematical Programming Computation},
	number = {3},
	pages = {219--260},
	publisher = {Springer},
	title = {Pyomo: modeling and solving mathematical programs in {P}ython},
	volume = {3},
	year = {2011}}

@article{duran1986outer,
  title={An outer-approximation algorithm for a class of mixed-integer nonlinear programs},
  author={Duran, Marco A and Grossmann, Ignacio E},
  journal={Mathematical Programming},
  volume={36},
  number={3},
  pages={307--339},
  year={1986},
  publisher={Springer}
}

@article{al1983jointly,
  title={Jointly constrained biconvex programming},
  author={Al-Khayyal, Faiz A and Falk, James E},
  journal={Mathematics of Operations Research},
  volume={8},
  number={2},
  pages={273--286},
  year={1983},
  publisher={INFORMS}
}

@Article{mccormick76,
  author        = {McCormick, Garth P.},
  title         = {Computability of Global Solutions to Factorable Nonconvex
                  Programs: Part {I} --- {C}onvex Underestimating Problems},
  journal       = {Mathematical Programming},
  year          = {1976},
  volume        = {10},
  number        = {1},
  pages         = {146-175}
}

@article{anstreicher2021convex,
  title={Convex hull representations for bounded products of variables},
  author={Anstreicher, Kurt M and Burer, Samuel and Park, Kyungchan},
  journal={Journal of Global Optimization},
  volume={80},
  number={4},
  pages={757--778},
  year={2021},
  publisher={Springer}
}

@article{nguyen2018deriving,
  title={Deriving convex hulls through lifting and projection},
  author={Nguyen, Trang T and Richard, Jean-Philippe P and Tawarmalani, Mohit},
  journal={Mathematical Programming},
  volume={169},
  number={2},
  pages={377--415},
  year={2018},
  publisher={Springer}
}

@article{belotti2010valid,
  title={Valid inequalities and convex hulls for multilinear functions},
  author={Belotti, Pietro and Miller, Andrew J and Namazifar, Mahdi},
  journal={Electronic Notes in Discrete Mathematics},
  volume={36},
  pages={805--812},
  year={2010},
  publisher={Elsevier}
}

@article{sherali2002enhancing,
  title={Enhancing {RLT} relaxations via a new class of semidefinite cuts},
  author={Sherali, Hanif D and Fraticelli, Barbara MP},
  journal={Journal of Global Optimization},
  volume={22},
  number={1},
  pages={233--261},
  year={2002},
  publisher={Springer}
}

@incollection{qualizza2011linear,
  title={Linear programming relaxations of quadratically constrained quadratic programs},
  author={Qualizza, Andrea and Belotti, Pietro and Margot, Fran{\c{c}}ois},
  booktitle={Mixed Integer Nonlinear Programming},
  pages={407--426},
  year={2011},
  publisher={Springer}
}

@book{sherali2013reformulation,
  title={A reformulation-linearization technique for solving discrete and continuous nonconvex problems},
  author={Sherali, Hanif D and Adams, Warren P},
  volume={31},
  year={2013},
  publisher={Springer Science \& Business Media}
}

@article{fortet1960applications,
  title={Applications de l'algebre de {B}oole en recherche op{\'e}rationelle},
  author={Fortet, Robert},
  journal={Revue Fran{\c{c}}aise de Recherche Op{\'e}rationelle},
  volume={4},
  number={14},
  pages={17--26},
  year={1960}
}

@article{vielma2017extended,
  title={Extended formulations in mixed integer conic quadratic programming},
  author={Vielma, Juan Pablo and Dunning, Iain and Huchette, Joey and Lubin, Miles},
  journal={Mathematical Programming Computation},
  volume={9},
  number={3},
  pages={369--418},
  year={2017},
  publisher={Springer}
}

@article{gally2018framework,
  title={A framework for solving mixed-integer semidefinite programs},
  author={Gally, Tristan and Pfetsch, Marc E and Ulbrich, Stefan},
  journal={Optimization Methods and Software},
  volume={33},
  number={3},
  pages={594--632},
  year={2018},
  publisher={Taylor \& Francis}
}

@article{burer2008finite,
  title={A finite branch-and-bound algorithm for nonconvex quadratic programming via semidefinite relaxations},
  author={Burer, Samuel and Vandenbussche, Dieter},
  journal={Mathematical Programming},
  volume={113},
  number={2},
  pages={259--282},
  year={2008},
  publisher={Springer}
}

@article{elloumi2019global,
  title={Global solution of non-convex quadratically constrained quadratic programs},
  author={Elloumi, Sourour and Lambert, Am{\'e}lie},
  journal={Optimization Methods and Software},
  volume={34},
  number={1},
  pages={98--114},
  year={2019},
  publisher={Taylor \& Francis}
}

@article{misener2013glomiqo,
  title={Glo{MIQO}: Global mixed-integer quadratic optimizer},
  author={Misener, Ruth and Floudas, Christodoulos A},
  journal={Journal of Global Optimization},
  volume={57},
  number={1},
  pages={3--50},
  year={2013},
  publisher={Springer}
}

@article{dong2018compact,
  title={Compact disjunctive approximations to nonconvex quadratically constrained programs},
  author={Dong, Hongbo and Luo, Yunqi},
  journal={arXiv preprint arXiv:1811.08122},
  year={2018}
}

@article{dong2016relaxing,
  title={Relaxing nonconvex quadratic functions by multiple adaptive diagonal perturbations},
  author={Dong, Hongbo},
  journal={SIAM Journal on Optimization},
  volume={26},
  number={3},
  pages={1962--1985},
  year={2016},
  publisher={SIAM}
}

@article{xia2020globally,
  title={Globally solving nonconvex quadratic programs via linear integer programming techniques},
  author={Xia, Wei and Vera, Juan C and Zuluaga, Luis F},
  journal={INFORMS Journal on Computing},
  volume={32},
  number={1},
  pages={40--56},
  year={2020},
  publisher={INFORMS}
}

@article{gondzio2021global,
  title={Global solutions of nonconvex standard quadratic programs via mixed integer linear programming reformulations},
  author={Gondzio, Jacek and Y{\i}ld{\i}r{\i}m, E Alper},
  journal={Journal of Global Optimization},
  volume={81},
  number={2},
  pages={293--321},
  year={2021},
  publisher={Springer}
}

@article{wei2024convex,
  title={On the convex hull of convex quadratic optimization problems with indicators},
  author={Wei, Linchuan and Atamt{\"u}rk, Alper and G{\'o}mez, Andr{\'e}s and K{\"u}{\c{c}}{\"u}kyavuz, Simge},
  journal={Mathematical Programming},
  volume={204},
  number={1},
  pages={703--737},
  year={2024},
  publisher={Springer}
}

@article{bienstock2025,
  title={Accurate linear cutting-plane relaxations for {ACOPF}},
  author={Bienstock, Daniel and Villagra, Mat{\'\i}as},
  journal={Mathematical Programming Computation},
  pages={79--133},
  volume={18},
  year={2026},
  publisher={Springer}
}

@article{brun2025,
  title={Alternating Methods for Large-Scale {AC} Optimal Power Flow with Unit Commitment},
  author={Brun, Matthew and Lee, Thomas and Lauinger, Dirk and Chen, Xin and Sun, Xu Andy},
  journal={arXiv preprint arXiv:2505.06211},
  year={2025}
}

@article{zhang2023,
  title={On solving unit commitment with alternating current optimal power flow on {GPU}},
  author={Zhang, Weiqi and Kim, Youngdae and Kim, Kibaek},
  journal={arXiv preprint arXiv:2310.13145},
  year={2023}
}

@techreport{elbert2024,
  title={{{ARPA-E} Grid Optimization ({GO}) Competition Challenge 3}},
  author={Elbert, Stephen and Holzer, Jesse and Veeramany, Arun and O'Neill, Richard and Mittelmann, Hans and Coffrin, Carleton and Garcia, Manuel and Parker, Robert and Elgindy, Tarek and Hale, Elaine and others},
  year={2024},
  institution={DOE Open Energy Data Initiative (OEDI)},
  address={Pacific Northwest National Laboratory},
  note= {\url{http://dx.doi.org/10.25984/2426334}}
}

@Article{GonzalezRodriguez2023,
  author    = {Brais Gonz{\'{a}}lez-Rodr{\'{\i}}guez and Joaqu{\'{\i}}n Ossorio-Castillo and Julio Gonz{\'{a}}lez-D{\'{\i}}az and {\'{A}}ngel M. Gonz{\'{a}}lez-Rueda and David R. Penas and Diego Rodr{\'{\i}}guez-Mart{\'{\i}}nez},
  journal   = {Journal of Global Optimization},
  title     = {Computational advances in polynomial optimization: {RAPOSa}, a freely available global solver},
  year      = {2023},
  number    = {3},
  pages     = {541--568},
  volume    = {85},
  doi       = {10.1007/s10898-022-01229-w},
  publisher = {Springer Science and Business Media {LLC}},
}

@article{raposaconic,
	author = {Brais Gonz{\'a}lez-Rodr{\'\i}guez and Ra{\'u}l Alvite-Paz{\'o} and Samuel Alvite-Paz{\'o} and Bissan Ghaddar and Julio Gonz{\'a}lez-D{\'\i}az},
	title = {Polynomial Optimization: Tightening {RLT}-Based Branch-and-Bound Schemes with Conic Constraints},
    journal = {Journal of Optimization Theory and Applications},
    volume = {204},
	year = {2025},
    number={1},
    pages={12}
}

@article{grone1984,
  title={Positive definite completions of partial {H}ermitian matrices},
  author={Grone, Robert and Johnson, Charles R and S{\'a}, Eduardo M and Wolkowicz, Henry},
  journal={Linear algebra and its applications},
  volume={58},
  pages={109--124},
  year={1984},
  publisher={Elsevier}
}

@manual{gurobi2025,
  title        = {Gurobi Optimizer 13.0.0 Reference Manual},
  author       = {{Gurobi Optimization, LLC}},
  year         = {2025},
  url          = {https://www.gurobi.com/documentation},
}

@manual{baron2025,
  title        = {BARON User Manual, Version 25.12.10},
  author       = {Sahinidis, N. V.},
  year         = {2025},
  url          = {https://www.minlp.com/baron-downloads},
}

@article{savelli2021,
  title={Electricity prices and tariffs to keep everyone happy: A framework for fixed and nodal prices coexistence in distribution grids with optimal tariffs for investment cost recovery},
  author={Savelli, Iacopo and Morstyn, Thomas},
  journal={Omega},
  volume={103},
  pages={102450},
  year={2021},
  publisher={Elsevier}
}

@article{amorosi2025,
  title={A time space network model for a truck and Multi-drone delivery system with battery recharging and variable speeds},
  author={Amorosi, Lavinia and Dell’Olmo, Paolo and Puerto, Justo and Valverde, Carlos},
  journal={Omega},
  pages={103399},
  year={2025},
  publisher={Elsevier}
}


\end{document}